# CRITICALLY LOADED QUEUEING MODELS THAT ARE THROUGHPUT SUBOPTIMAL[1]

By Rami Atar and Gennady Shaikhet

*Technion—Israel Institute of Technology and Carnegie Mellon University*

This paper introduces and analyzes the notion of throughput suboptimality for many-server queueing systems in heavy traffic. The queueing model under consideration has multiple customer classes, indexed by a finite set $\mathcal{I}$, and heterogenous, exponential servers. Servers are dynamically chosen to serve customers, and buffers are available for customers waiting to be served. The arrival rates and the number of servers are scaled up in such a way that the processes representing the number of class-$i$ customers in the system, $i \in \mathcal{I}$, fluctuate about a static fluid model, that is assumed to be critically loaded in a standard sense. At the same time, the fluid model is assumed to be throughput suboptimal. Roughly, this means that the servers can be allocated so as to achieve a total processing rate that is greater than the total arrival rate. We show that there exists a dynamic control policy for the queueing model that is efficient in the following strong sense: Under this policy, for every finite $T$, the measure of the set of times prior to $T$, at which at least one customer is in the buffer, converges to zero in probability as the arrival rates and number of servers go to infinity. On the way to prove our main result, we provide a characterization of throughput suboptimality in terms of properties of the buffer-station graph.

**1. Introduction.** In this paper, we study a class of many-server queueing systems in heavy traffic, that are critically loaded in a standard sense, but exhibit a behavior that is typical to subcritically loaded systems. We introduce the notion of throughput suboptimality for an underlying fluid model, and show that it plays a central role in determining and explaining this behavior.

Received February 2007; revised April 2008.

[1]Supported in part by the Israel Science Foundation (Grant no. 126/02) and the Technion fund for the promotion of research.

*AMS 2000 subject classifications.* 60K25, 68M20, 90B22, 90B36, 60F05.

*Key words and phrases.* Multi-class queueing systems, heavy traffic, scheduling and routing, throughput optimality, asymptotic null controllability, buffer-station graph.

The queueing model under consideration has multiple customer classes, indexed by a finite set $\mathcal{I}$, and heterogeneous exponential servers. The servers are grouped in pools, indexed by a finite set $\mathcal{J}$, and it is assumed that each pool has a large number of servers that have identical capabilities. In particular, the rate, denoted by $\mu_{ij}$, at which a server from pool $j \in \mathcal{J}$ serves customers from class $i \in \mathcal{I}$, depends on both $i$ and $j$. It is also possible that servers from pool $j$ cannot serve class-$i$ customers, in which case we write $\mu_{ij} = 0$. The arrival of class-$i$ customers is modeled as a renewal process with rate $\lambda_i$, $i \in \mathcal{I}$. Servers are dynamically chosen to serve customers, and buffers are available to accommodate customers that wait to be served (see Figure 1). The model is considered in a many-server heavy traffic regime, in which the number of servers at each pool and the arrival rates are scaled up at a nearly fixed proportion, and in such a way that the processes that represent the number of class-$i$ customers in the system, $i \in \mathcal{I}$, fluctuate about a certain static fluid model. This fluid model is assumed to be critically loaded, in a standard sense. In particular, (1) servers can be allocated in such a way that the total processing rate devoted to class-$i$ customers is equal to the arrival rate $\lambda_i$, for every $i \in \mathcal{I}$; and (2) property (1) does not hold if one of the arrival rates $\lambda_i$ is replaced by some $\lambda_i' > \lambda_i$ (there are some further assumptions; see Section 2). It is possible for such a model to satisfy the following condition: servers can be allocated so as to achieve a total processing rate that is greater than the total arrival rate (see Section 2 for a precise statement). If this condition holds, we say that the fluid model is *throughput suboptimal*. Our main result (Theorem 1) shows that when the fluid model is throughput suboptimal, one can find a dynamic control policy for the queueing model that exhibits a strong form of efficiency: Under this policy, for every finite $T$, the measure of the set of times prior to $T$, at which at least one customer is in the buffer, converges to zero in probability as the arrival rates and number of servers go to infinity. Thus, although the system is critically loaded, its buffers are "essentially" empty, as if the system is subcritically loaded (see [5], Theorem 6.8(i), for a typical asymptotic result in a subcritical regime in diffusion scale, where the buffers are empty in the limit).

This work is motivated by recent progress on scheduling and routing control problems for many-server systems, and their diffusion-scale limits; see references cited in [1], Section 2.3.3. In many of these references, one attempts to find a dynamic routing policy that minimizes a given performance measure. Viewing these hard problems at a diffusion scale typically simplifies the task, because at the scaling limit the problem transforms into that of optimally controlling a diffusion process. Diffusion scale is natural because at this scale the primitive processes (associated with arrival and service) exhibit nontrivial random fluctuations. To explain how our result is related to this viewpoint, let $n$ denote a parameter by which the number of servers and the arrival rates are scaled [as in (2.7); see Section 2 for the precise setting]. Denote the total queue-length of the $n$th system by $Q^n$ (starting in Section 2 this process will be denoted by $e \cdot Y^n$). A typical control problem is that of minimizing $E[\int_0^T \widehat{Q}^n(t)\,dt]$, where $\widehat{Q}^n = n^{-1/2}Q^n$, a diffusion-scale version of $Q^n$. Our result can be viewed as a contribution to this line of work. Indeed, assuming throughput suboptimality of the fluid model, it implies that there exists a policy under which, for every $T < \infty$, the empirical law associated with $\widehat{Q}^n|_{[0,T]}$ converges weakly to the unit point mass at zero, as $n \to \infty$ (although it does not address an expected cost of the above type). However, we do not present the result in diffusion scale, because it is sharper and, in fact, establishes the above statement about the empirical law of the process $Q^n|_{[0,T]}$ itself.

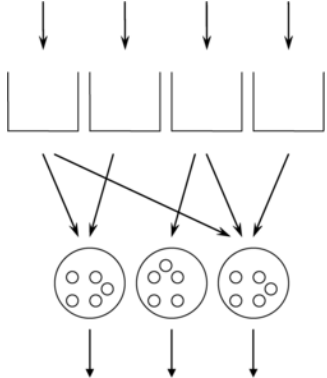

Fig. 1. *A queueing model with four customer classes and three service pools.*

A related analysis appears in [3], where the same model is proved to satisfy a stronger result under a stronger assumption. While the current paper addresses the capability to maintain a system with no customer in the buffer "most of the time," with large probability, the result of [3] concerns maintaining a system with no customers in the buffer "at all times" (apart from an initial transition phase), with large probability. More precisely, under appropriate assumptions, it is shown in [3] that there exists a policy under which, for every $0 < \varepsilon < T < \infty$, the probability that at least one customer is present in the buffer any time within $[\varepsilon, T]$ approaches 0 in the scaling limit. This phenomenon is shown to be related, on one hand to a formulation of the limiting diffusion model as a diffusion with singular control [3], Section 3. On the other hand, it is shown to be related to a condition on the graph that encodes the network's structure. This graph has a vertex for each class $i \in \mathcal{I}$, a vertex for each server pool $j \in \mathcal{J}$, and an edge, with an associated weight $\mu_{ij}$, between a class vertex $i$ and a pool vertex $j$ if, and only if $\mu_{ij} > 0$. The assumption of [3] is the existence of a cycle $p$ in this graph, having a negative total signed weight, $\mu(p)$, where the ($p$-dependent) signs of the weights $\mu_{ij}$ are appropriately defined [as in equation (3.1) in the current paper; see also equation (3.2) for a definition of $\mu(p)$ as the sum of the signed weights along $p$]. We will show that the algebraic condition alluded to above is a special case of the main assumption of the current paper, namely throughput suboptimality. We will also characterize the latter condition in terms of the graph and the signed weights, and show that throughput suboptimality may occur in one of two ways: The existence of either a cycle or an open path $p$ (appropriately defined), with signed weight $\mu(p) < 0$ (Theorem 2).

We make two further remarks about the relation to [3]. First, the difference between having no customers in the buffer for a given period of time $[\varepsilon, T]$ (as in [3]) and having no customers in the buffer most of the time within $[\varepsilon, T]$, may be significant with regard to the queue-length performance measure. In fact, under the policy constructed in the current paper, there are short time periods in which large queues build. We believe that a result of the type of [3] is not possible under the conditions of the current paper, but we do not prove this claim. Second, the results of [3] allow for both preemptive policies (where service to a customer can be interrupted and resumed at a later time, possibly at a different server) and nonpreemptive ones (where service cannot be interrupted), while the current paper only treats preemptive policies. We leave open the question of whether analogous results are possible for the nonpreemptive case.

The results of [3] and the current paper reveal two aspects of a phenomenon, where critically loaded many-server systems behave as subcritically loaded. As our main result shows, the notion of throughput suboptimality captures this phenomenon. It is reasonable to expect that this connection continues to hold in a wider family of critically loaded many server models, with more general structure.

The main tool in analyzing the probabilistic model is a related deterministic dynamic fluid model, that roughly replaces stochastic fluctuations by deterministic ones. Throughput suboptimality is shown to have an effect on this model that is similar to the one discussed above, where quantities that represent queue-lengths are shown to be small, in an appropriate sense (see Theorem 3). The proof of the result relies on the graph-theoretic characterization alluded to above, and specifically uses the existence of a path $p$ with the property $\mu(p) < 0$. The result for the probabilistic model follows from the deterministic one in a relatively straightforward way.

The organization of the paper is as follows. Section 2 contains the description of the model and assumptions, and the statement of the main result. Some numerical examples are given at the end of this section. Section 3 provides an algebraic characterization of throughput (sub) optimality. The dynamic fluid model is introduced in Section 4. A property for this model that is analogous to the main result is proved, based on the results of Section 3. Relying on the deterministic model results, we provide in Section 5 a proof of the main result.

NOTATION. Write $\mathbb{N}$ for the set of positive integers, $\mathbb{Z}_+$ for the set of nonnegative integers, and $\mathbb{R}_+ = [0, \infty)$. For a real number $a$, let $a^+ = \max\{0, a\}$, $a^- = -\min\{0, a\}$. For a positive integer $d$ and $x \in \mathbb{R}^d$, let $\|x\| = \sum_{i=1}^d |x_i|$. For $v, u \in \mathbb{R}^d$ let $v \cdot u = \sum_{i=1}^d u_i v_i$. The symbols $e_i$ denote the unit coordinate vectors and $e = (1, \ldots, 1)$. The dimension of $e$ may change from one expression to another. Denote by $\mathbb{D}(\mathbb{R}^d)$ the space of all cadlag functions (i.e., right continuous and having finite left limits) from $\mathbb{R}_+$ to $\mathbb{R}^d$. Endow $\mathbb{D}(\mathbb{R}^d)$ with the usual Skorohod topology (cf. [4]). If $X^n$, $n \in \mathbb{N}$, and $X$ are processes with sample paths in $\mathbb{D}(\mathbb{R}^d)$, write $X^n \Rightarrow X$ to denote weak convergence of the measures induced by $X^n$ [on $\mathbb{D}(\mathbb{R}^d)$] to the measure induced by $X$. Denote $|X|_t^* = \sup_{0 \leq u \leq t} |X(u)|$ for $X \in \mathbb{D}(\mathbb{R})$ and $\|X\|_t^* = \sup_{0 \leq u \leq t} \|X(u)\|$ for $X \in \mathbb{D}(\mathbb{R}^d)$.

## 2. Setting and main result.

2.1. *Probabilistic queueing model.* A precise description of the queueing model is as follows. A complete probability space $(\Omega, \mathcal{F}, \mathbb{P})$ is given, supporting all stochastic processes defined below. Expectation with respect to $\mathbb{P}$ is denoted by $\mathbb{E}$. The queueing model is parameterized by $n \in \mathbb{N}$. It has $I$ customer classes and $J$ service stations. Station $j$ has $N_j^n$ identical servers. The classes are labeled as $1, \ldots, I$ and the stations as $I + 1, \ldots, I + J$:

$$\mathcal{I} = \{1, \ldots, I\}, \qquad \mathcal{J} = \{I + 1, \ldots, I + J\}.$$

Arrivals are modeled as renewal processes with finite second moment for the interarrival time. More precisely, we are given parameters $\lambda_i^n > 0$, $i \in \mathcal{I}$, $n \in \mathbb{N}$, and independent sequences of strictly positive i.i.d. random variables $\{\breve{U}_i(k), k \in \mathbb{N}\}$, $i \in \mathcal{I}$, with mean $E\breve{U}_i(1) = 1$ and squared coefficient of variation $C_{U,i}^2 = (E\breve{U}_i(1))^{-2} \operatorname{Var}(\breve{U}_i(1)) \in [0, \infty)$. With $\sum_1^0 = 0$, the number of class-$i$ arrivals up to time $t$ at the $n$th system is given as

$$A_i^n(t) = \sup\left\{l \geq 0 : \sum_{k=1}^l \frac{\breve{U}_i(k)}{\lambda_i^n} \leq t\right\}, \qquad t \geq 0.$$

For $i \in \mathcal{I}, j \in \mathcal{J}$ and $n \in \mathbb{N}$ we are given parameters $\mu_{ij}^n \geq 0$, representing the service rate of a class-$i$ customer by a server of station $j$. There is a possibility for $\mu_{ij}^n = 0$, in which case we say that class-$i$ customers cannot be served at station $j$. For every $(i, j) \in \mathcal{I} \times \mathcal{J}$, we denote by $\Psi_{ij}^n(t)$ the number of class-$i$ customers being served in station $j$ at time $t$. By definition,

$$(2.1) \qquad \Psi_{ij}^n(t) = 0 \qquad \text{for } (i, j) \text{ s.t. } \mu_{ij}^n = 0.$$

Service times are modeled as independent exponential random variables. To this end, let $S_{ij}^n$, $(i,j) \in \mathcal{I} \times \mathcal{J}$, be Poisson processes with rate $\mu_{ij}^n$ (where a Poisson process of zero rate is the zero process), mutually independent and independent of the arrival processes. Let $T_{ij}^n(t)$ denote the time up to $t$ devoted to a class-$i$ customer by a server, summed over all servers from station $j$, and note that

$$T_{ij}^n(t) = \int_0^t \Psi_{ij}^n(s)\,ds, \qquad i \in \mathcal{I}, j \in \mathcal{J}, t \geq 0.$$

The number of service completions of class-$i$ customers by all servers of station $j$ by time $t$ is, by assumption, given by $D_{ij}^n(t) := S_{ij}^n(T_{ij}^n(t))$. See [3] for explanation on the exponential service time property of this model. We refer to $D^n = (D_{ij}^n, (i,j) \in \mathcal{I} \times \mathcal{J})$ as the *departure process*. The processes $A^n$ and $S^n$ will be referred to as the *primitive processes*.

Denoting by $X_i^n(t)$ the number of class-$i$ customers in the system at time $t$, and setting $X_i^{0,n} = X_i^n(0)$, it is clear from the above that

$$(2.2) \quad X_i^n(t) = X_i^{0,n} + A_i^n(t) - \sum_{j \in \mathcal{J}} S_{ij}^n\left(\int_0^t \Psi_{ij}^n(s)\,ds\right), \qquad i \in \mathcal{I}, t \geq 0.$$

For simplicity, the initial conditions $X_i^{0,n}$ are assumed to be deterministic. Finally, we introduce the processes $Y_i^n(t)$, representing the number of class-$i$ customers that are in the queue (and not being served) at time $t$, and $Z_j^n(t)$, representing the number of servers at station $j$ that are idle at time $t$. Clearly, we have the following relations:

$$(2.3) \qquad\qquad Y_i^n(t) + \sum_{j \in \mathcal{J}} \Psi_{ij}^n(t) = X_i^n(t), \qquad i \in \mathcal{I},$$

$$(2.4) \qquad\qquad Z_j^n(t) + \sum_{i \in \mathcal{I}} \Psi_{ij}^n(t) = N_j^n, \qquad j \in \mathcal{J}.$$

Also, the following holds by definition:

$$(2.5) \quad Y_i^n(t) \geq 0, \qquad Z_j^n(t) \geq 0, \qquad \Psi_{ij}^n(t) \geq 0, \qquad i \in \mathcal{I}, j \in \mathcal{J}, t \geq 0.$$

Equations (2.1)–(2.5) indicate some properties of the processes involved, but they do not characterize these processes, because the process $\Psi^n$ has not yet been described. As reflected in the following definition, we regard $\Psi^n$ as a control process that can be obtained as "feedback" from the "state" process $X^n$ and the arrival process $A^n$.

DEFINITION 1. Fix $n$. We say that a process $\Psi^n = (\Psi_{i,j}^n)_{j \in \mathcal{J}}^{i \in \mathcal{I}}$ where, for $(i,j) \in \mathcal{I} \times \mathcal{J}$, $\Psi_{i,j}^n$ takes values in $\mathbb{Z}_+$ and has right-continuous paths, is a *scheduling control policy* (SCP) if the following conditions hold:

(i) Given initial data $X^{n,0}$ and primitive processes $A^n$ and $S^n$, there exist processes $X^n$, $Y^n$ and $Z^n$ with values in $\mathbb{Z}_+^{\mathcal{I}}$, $\mathbb{Z}_+^{\mathcal{I}}$ and $\mathbb{Z}_+^{\mathcal{J}}$, respectively, such that (2.1)–(2.5) are met;

(ii) For every $t \geq 0$, $\Psi^n(t)$ is measurable on $\sigma\{X^n(s), A^n(s) : s \leq t\}$.

Note that uniqueness of the processes $X^n$, $Y^n$ and $Z^n$, given $A^n$, $S^n$ and $\Psi^n$, is immediate from (2.2)–(2.4). Note also that according to this definition, service to a customer can be stopped and resumed at a later time, possibly in a different station.

We will use some elementary graph theoretic terminology and notation as follows (see, e.g., [6] for standard definitions). For a nonempty set $V$ and $E \subseteq V \times V$, we write $G = (V, E)$ for the graph with vertex set $V$ and edge set $E$. A vertex having exactly one neighbor is called a *leaf vertex*, and an edge joining a leaf vertex is called a *leaf edge*. A connected graph that does not contain cycles is called a *tree*.

We denote the index set for all customer classes and service stations by $\mathcal{V} := \mathcal{I} \cup \mathcal{J}$, and the set of all class–station pairs by $\mathcal{E} := \mathcal{I} \times \mathcal{J}$. Some of the elements $(i, j)$ of $\mathcal{E}$ correspond to class–station pairs where station $j$ can serve class $i$. We encode this information by setting

$$(2.6) \qquad \mathcal{E}_{\mathrm{a}} = \{(i,j) \in \mathcal{I} \times \mathcal{J} : \mu_{ij}^n > 0\},$$

where throughout we assume that $\mathcal{E}_{\mathrm{a}}$ does not depend on $n$. A class–station pair $(i, j) \in \mathcal{E}_{\mathrm{a}}$ is said to be an *activity*. Throughout, if $\mathcal{E}_1$ is a subset of $\mathcal{E}$, we write $\mathcal{E}_1^c$ for the complement of $\mathcal{E}_1$ with respect to $\mathcal{E}$. The set of class–station pairs that are not activities is denoted by $\mathcal{E}_{\mathrm{a}}^c \equiv \mathcal{E} \setminus \mathcal{E}_{\mathrm{a}}$. We denote $\mathcal{G}_{\mathrm{a}} = (\mathcal{V}, \mathcal{E}_{\mathrm{a}})$, and refer to it as the graph of activities.

### 2.2. *Static fluid model: heavy traffic and throughput optimality.*

*Heavy traffic condition and related assumptions.* We will assume that the parameters of the probabilistic queueing model satisfy certain conditions that indicate that the system is critically loaded. To specify these conditions, we will introduce a deterministic, static fluid model, defined in terms of a simple linear program.

We assume that there are constants $\lambda_i, \nu_j \in (0, \infty)$, $i \in \mathcal{I}$, $j \in \mathcal{J}$, and $\mu_{ij} \in (0, \infty)$, $(i, j) \in \mathcal{E}_{\mathrm{a}}$, such that

$$(2.7) \qquad \begin{aligned} n^{-1}\lambda_i^n &\to \lambda_i, & i &\in \mathcal{I}, \\ n^{-1}N_j^n &\to \nu_j, & j &\in \mathcal{J}, \\ \mu_{ij}^n &\to \mu_{ij}, & (i,j) &\in \mathcal{E}_{\mathrm{a}}. \end{aligned}$$

We set $\mu_{ij} = 0$ for $(i, j) \in \mathcal{E}_{\mathrm{a}}^c$. [Later on we strengthen (2.7) above; cf. Assumption 3.]

Consider a fluid model, where the arrival and service processes are replaced by deterministic flows with corresponding rates $\lambda_i$ and $\mu_{ij}$. There are $I$ classes of incoming fluid and $J$ processing stations (while for the probabilistic queueing model we used the terms *class-i customers* and *service station j*, for the current model we use *class-i fluid* and *processing station j*). Station $j$ has capacity to hold $\nu_j$ units of fluid. When station $j$ contains $\psi_{ij}$ units of class-$i$ fluid, for all $i \in \mathcal{I}$ and $j \in \mathcal{J}$, the rate at which class-$i$ fluid is processed at station $j$ (and leaves the station) is $\mu_{ij}\psi_{ij}$. The overall rate at which class-$i$ fluid is processed (and leaves the system) is $\sum_j \mu_{ij}\psi_{ij}$. Let $\Xi$ be the set of $I \times J$ matrices $\xi$ with $\xi_{ij} \geq 0$, $(i,j) \in \mathcal{E}$, and $\sum_i \xi_{ij} \leq 1$, $j \in \mathcal{J}$. For $\xi \in \Xi$, $\xi_{ij}$ will represent the fraction of the service capacity from station $j$ allocated to process class-$i$ fluid. We call an element of $\Xi$ an *allocation matrix*. The fluid model uses a fixed allocation matrix for all times (hence the term "static" model). Set $\bar{\mu}_{ij} = \mu_{ij}\nu_j$, $(i,j) \in \mathcal{E}$. Note that $\xi_{ij}$ is equal to the amount of class-$i$ fluid contained in station $j$, normalized by the capacity of station $j$, namely $\psi_{ij}/\nu_j$. Consequently, $\mu_{ij}\psi_{ij} = \bar{\mu}_{ij}\xi_{ij}$. Consider the following linear program:

*Find $\{\xi_{ij}, (i,j) \in \mathcal{E}\}$ and $\rho \in \mathbb{R}_+$ so as to minimize $\rho$ subject to*

(2.8)
$$\begin{cases} \sum_{j \in \mathcal{J}} \bar{\mu}_{ij}\xi_{ij} = \lambda_i, & i \in \mathcal{I}, \\ \sum_{i \in \mathcal{I}} \xi_{ij} \leq \rho, & j \in \mathcal{J}, \\ \xi_{ij} \geq 0, & (i,j) \in \mathcal{E}. \end{cases}$$

For $\rho \in [0,1]$, a $\xi$ as above is clearly an allocation matrix. The first line of (2.8) expresses that the system is balanced, in the sense that for each $i$, the total processing rate of class-$i$ fluid equals the rate at which fluid of this class enters the system. We will assume throughout that the system is critically loaded, in the following sense.

ASSUMPTION 1. There exists an optimal solution $(\xi^*, \rho^*)$ to the linear program (2.8), satisfying $\sum_{i \in \mathcal{I}} \xi_{ij}^* = 1$ for all $j \in \mathcal{J}$ (and consequently $\rho^* = 1$).

This assumption is weaker than the *heavy traffic condition* of [7] where the optimal solution is assumed, in addition, to be unique. Throughout the paper, we fix $\xi_{ij}^*$ satisfying Assumption 1, and also let

(2.9)           $\psi_{ij}^* = \xi_{ij}^*\nu_j, \qquad x_i^* = \sum_j \xi_{ij}^*\nu_j, \qquad i \in \mathcal{I}, j \in \mathcal{J}.$

The following simple relations follow directly from the above assumption:

(2.10)      $\sum_{j \in \mathcal{J}} \psi_{ij}^* = x_i^*, \qquad \sum_{i \in \mathcal{I}} \psi_{ij}^* = \nu_j, \qquad \lambda_i = \sum_{j \in \mathcal{J}} \mu_{ij}\psi_{ij}^*, \qquad i \in \mathcal{I}, j \in \mathcal{J}.$

The quantity $\psi_{ij}^*$ represents the amount of class-$i$ fluid that station $j$ contains (and processes), under the allocation matrix $\xi^*$, and $x_i^*$ represents the total amount of class-$i$ fluid being processed.

Following [7], an activity $(i, j) \in \mathcal{E}_a$ is said to be *basic* (resp., *nonbasic*) if $\xi_{ij}^* > 0$ (resp., $= 0$). Define the graph of basic activities $\mathcal{G}_{ba}$ to be the subgraph of $\mathcal{G}_a$ having $\mathcal{V}$ as a vertex set, and the collection

$$\mathcal{E}_{ba} := \{(i, j) \in \mathcal{E}_a : \xi_{ij}^* > 0\}$$

of basic activities as an edge set. The following will be assumed throughout.

ASSUMPTION 2. *The graph $\mathcal{G}_{ba}$ is a tree.*

Under the heavy traffic condition of [7], alluded to above, Assumption 2 above is known [8] to be equivalent to the *complete resource pooling condition* [3, 7], of which one of the equivalent formulations is that all vertices in $\mathcal{J}$ communicate via edges in $\mathcal{G}_{ba}$, expressing a strong mode of cooperation between the service stations. Note that the combination of Assumptions 1 and 2 forms a weaker assumption than the combination of the heavy traffic and complete resource pooling conditions (assumed, e.g., in [3]). It is possible for the linear program to have multiple optimal solutions, and it is only required that one of them satisfies the tree structure.

*Throughput optimality.* Assumption 1 expresses critical load on the system, in the sense described earlier, but it does not exclude the possibility that the total processing rate can exceed the total arrival rate. Namely, it is possible that there exists an allocation matrix $\xi$ under which

$$(2.11) \qquad \sum_{(i,j) \in \mathcal{E}} \bar{\mu}_{ij} \xi_{ij} > \sum_{i \in \mathcal{I}} \lambda_i.$$

The set of allocation matrices $\xi \in \Xi$ that satisfy

$$(2.12) \qquad \sum_{j \in \mathcal{J}} \xi_{ij} \nu_j \le x_i^* \qquad \text{for all } i \in \mathcal{I}$$

is of interest. Under these allocation matrices, for each $i \in \mathcal{I}$, the total amount of class-$i$ fluid being processed does not exceed that under $\xi^*$. A condition involving simultaneously (2.11) and (2.12) will be key in this paper. We will say that the static fluid model is *throughput optimal* if the following holds:

$$(2.13) \qquad \begin{aligned} &\textit{Whenever } \xi \in \Xi \textit{ and } \sum_{j \in \mathcal{J}} \xi_{ij} \nu_j \le x_i^* \; \forall i \in \mathcal{I}, \\ &\textit{one has } \sum_{(i,j) \in \mathcal{E}} \bar{\mu}_{ij} \xi_{ij} \le \sum_{i \in \mathcal{I}} \lambda_i. \end{aligned}$$

We will say that the static fluid model is *throughput suboptimal* if it is not throughput optimal.

When the static fluid model is throughput suboptimal, one can find $\xi \in \Xi$ meeting (2.11) and (2.12). Recall that $x^*$ represents the amount of fluid of each class being processed in all stations. Thus, when (2.13) fails to hold, one can keep the same amount of fluid of each class as under $\xi^*$, and redistribute it among the stations so as to obtain a greater total processing rate than under $\xi^*$. Note, however, there is no guarantee that under $\xi$ the processing rate is sufficient for handling arrivals of all classes. Namely, it is possible that

$$(2.14) \qquad \text{there is a class } i \in \mathcal{I} \text{ for which} \qquad \sum_{j \in \mathcal{J}} \bar{\mu}_{ij} \xi_{ij} < \lambda_i,$$

and thus a use of $\xi$ may result in instability. [It is not hard to see that in the case where the optimal solution to the the linear program is unique, (2.14) holds for any $\xi \neq \xi^*$.] In the probabilistic queueing model, however, one can vary the allocation over time, and the existence of $\xi$ as above turns out to have a crucial impact. This is expressed in Theorem 1 below, which is our main result.

The following assumption regards the second order behavior of the parameters and initial condition.

ASSUMPTION 3. *There is a constant $c < \infty$ such that for all $i \in \mathcal{I}$, $j \in \mathcal{J}$ and $n \in \mathbb{N}$,*

$$(2.15) \quad |n^{-1}\lambda_i^n - \lambda_i| \vee |\mu_{ij}^n - \mu_{ij}| \vee |n^{-1}N_j^n - \nu_j| \vee |n^{-1}X_i^{0,n} - x_i^*| \leq cn^{-1/2}.$$

THEOREM 1. *Let Assumptions 1–3 hold. If the static fluid model is throughput suboptimal, then there exists a sequence of SCPs, under which for any fixed $0 < T < \infty$ and $\varrho > 1/2$,*

$$(2.16) \qquad \int_0^T 1_{\{e \cdot Y^n(s) > 0\}}\, ds \to 0 \qquad \text{in probability, as } n \to \infty,$$

$$(2.17) \qquad n^{-\varrho}\|X^n - X^{0,n}\|_T^* \to 0 \qquad \text{in probability, as } n \to \infty.$$

REMARK 1. Assumption 3 could be somewhat relaxed by replacing the bound on the last term in (2.15), namely $|n^{-1}X_i^{0,n} - x_i^*| \leq cn^{-1/2}$, by $(n^{-1}X_i^{0,n} - x_i^*) \leq cn^{-1/2}$, $n \in \mathbb{N}$, $i \in \mathcal{I}$, so as to cover cases where the initial load on the system is subcritical. There is a simple argument by which this can be deduced from Theorem 1, where one introduces virtual customers at time zero, thus increasing the value of $X^{0,n}$ so that (2.15) holds, and Theorem 1 is in force. We leave out the details.

The main idea of how throughput suboptimality is used in the proof of Theorem 1 is as follows. Recall that when the static fluid model is throughput suboptimal, one can find an allocation matrix $\xi \in \Xi$ under which the total processing rate is strictly greater than the total arrival rate $\sum_i \lambda_i$. In the probabilistic queueing model, we can interpret this condition by saying that when, for each $j$, the proportion of pool-$j$ servers allocated to class $i$ is (roughly) $\xi_{ij}$, the system operates with a total service rate that exceeds the total arrival rate. Since (2.14) may hold, a constant use of such proportions for a long time might result in instability, in the sense that at least one of the queues must build up. Therefore, such a strategy will not achieve (2.16). A slightly more sophisticated strategy is one that alternates between policies (i) and (ii), as follows. (i) Use, as above, proportions according to $\xi$, until the system reaches a state where the total number of customers is less than the total number of servers, but possibly with some queues building up; (ii) Rearrange customers in the system so that all buffers are empty. This is made possible because the number of customers is less than the total number of servers, and the vector of normalized numbers of customers is close to $x^*$ [note also that (2.12) is necessary because in the rearrangement process (ii), the total number of servers allocated to serve class $i$ is limited by the number class-$i$ customers in the system]. Then use a control that merely keeps the buffers empty for some time. The proof shows that this can be done in such a way that the length of the time intervals where (i) is applied are much shorter than those where (ii) is applied, so that (2.16) is achieved in the limit.

In order for our argument to be valid, there must exist an allocation matrix $\xi$ which achieves the inequality that is a part of the definition of throughput suboptimality, namely that under $\xi$, the total processing rate exceeds the total arrival rate. That, of course, does not imply that throughput suboptimality is necessary for the validity of (2.16), and so the question of whether a converse to our main result holds, is left open. We intend to address this question in future work.

### 2.3. *Examples.*

We demonstrate throughput suboptimality by some numerical examples.

EXAMPLE 1. Consider the following static fluid model in heavy traffic, with 2 classes of customers and 3 stations

$$\nu = \begin{pmatrix} 1 \\ 1 \\ 1 \end{pmatrix}, \qquad \lambda = \begin{pmatrix} 8 \\ 4 \end{pmatrix}, \qquad \mu = \bar{\mu} = \begin{pmatrix} 3 & 10 & 1 \\ 1 & 4 & 2 \end{pmatrix}.$$

The resulting optimal static allocation is as follows (2.9):

$$\psi^* = \xi^* = \begin{pmatrix} 1 & 0.5 & 0 \\ 0 & 0.5 & 1 \end{pmatrix} \quad \text{and} \quad x^* = \begin{pmatrix} 1.5 \\ 1.5 \end{pmatrix}.$$

To see that the static fluid model is throughput suboptimal, let $\varepsilon > 0$ be sufficiently small and consider the allocation matrix

$$\widehat{\xi} = \begin{pmatrix} 1-\varepsilon & 0.5+\varepsilon & 0 \\ \varepsilon & 0.5-\varepsilon & 1 \end{pmatrix}.$$

Clearly, we have $\sum_j \widehat{\xi}_{ij} \nu_j = x_i^*$ for every $i$. However, $\sum_{(i,j) \in \mathcal{E}} \widehat{\xi}_{ij} \bar{\mu}_{ij} > \lambda_1 + \lambda_2$. Thus, the condition of throughput optimality (2.13) is not satisfied. The result of Theorem 1 holds. We note that the assumptions of [3] are also valid in this example. See more information on this example at the end of Section 3.

EXAMPLE 2. In this example, the data is the same as in Example 1 above, except for one entry:

$$\nu = \begin{pmatrix} 1 \\ 1 \\ 1 \end{pmatrix}, \qquad \lambda = \begin{pmatrix} 8 \\ 4 \end{pmatrix}, \qquad \mu = \bar{\mu} = \begin{pmatrix} 3 & 10 & 1 \\ 0 & 4 & 2 \end{pmatrix}.$$

The resulting optimal static allocation is as follows:

$$\psi^* = \xi^* = \begin{pmatrix} 1 & 0.5 & 0 \\ 0 & 0.5 & 1 \end{pmatrix}, \qquad x^* = \begin{pmatrix} 1.5 \\ 1.5 \end{pmatrix}.$$

With $\varepsilon > 0$ sufficiently small, the matrix

$$\widehat{\xi} = \begin{pmatrix} 1-\varepsilon & 0.5+\varepsilon & 0 \\ 0 & 0.5-\varepsilon & 1 \end{pmatrix}$$

is an allocation matrix. Moreover, $\sum_j \widehat{\xi}_{ij} \nu_j = x_i^*$ and $\sum_{(i,j) \in \mathcal{E}} \widehat{\xi}_{ij} \bar{\mu}_{ij} > \lambda_1 + \lambda_2$. Thus, the static fluid model is throughput suboptimal. As shown in the end of Section 3, the conditions of [3] are not satisfied for this example.

EXAMPLE 3. Consider

$$\nu = \begin{pmatrix} 1 \\ 1 \\ 1 \end{pmatrix}, \qquad \lambda = \begin{pmatrix} 4 \\ 1 \\ 2 \end{pmatrix}, \qquad \mu = \bar{\mu} = \begin{pmatrix} 2 & 4 & 0.5 \\ 0.3 & 1 & 1 \\ 0.1 & 0.5 & 4 \end{pmatrix}.$$

The resulting optimal static allocation is as follows:

$$\psi^* = \xi^* = \begin{pmatrix} 1 & 0.5 & 0 \\ 0 & 0.5 & 0.5 \\ 0 & 0 & 0.5 \end{pmatrix} \quad \text{and} \quad x^* = \begin{pmatrix} 1.5 \\ 1 \\ 0.5 \end{pmatrix}.$$

The static fluid model for this example is throughput optimal, as we show at the end of Section 3, using the tools we develop in Section 3.

**3. Characterization of throughput optimality.** The main result of this section (Theorem 2) characterizes throughput optimality in terms of some graph-theoretic properties of the network. To state it, we need some definitions.

Recall that by Assumption 2, the graph $\mathcal{G}_{\mathrm{ba}}$ is a tree, and because by construction, it is a subgraph of $\mathcal{G}_{\mathrm{a}}$, all its edges are of the form $(i, j)$ where $i \in \mathcal{I}$ and $j \in \mathcal{J}$. In the definition below and elsewhere in this section, it will be convenient to identify $(i, j)$ with $(j, i)$ (where $i \in \mathcal{I}$ and $j \in \mathcal{J}$) when referring to an element of the edge set $\mathcal{E}$. Although the notation is abused, there will be no confusion since $\mathcal{I}$ and $\mathcal{J}$ do not intersect.

DEFINITION 2. (i) A subgraph $q = (\mathcal{V}_q, \mathcal{E}_q)$ of $\mathcal{G}_{\mathrm{ba}}$ is called a *basic path* if one has $\mathcal{V}_q = \{i_0, j_0, \ldots, i_k, j_k\}$, and

$$\mathcal{E}_q = \{(i_0, j_0), (j_0, i_1), \ldots, (i_k, j_k)\}$$

where $k \geq 1$ and $i_0, \ldots, i_k \in \mathcal{I}$, $j_0, \ldots, j_k \in \mathcal{J}$ are $2k + 2$ distinct vertices. Note that every edge of a basic path is a basic activity (i.e., an element of $\mathcal{E}_{\mathrm{ba}}$). Denote by $BP$ the set of basic paths. Basic paths are used in this paper mainly in order to define simple paths, as follows.

(ii) Let the leaves $i_0$ and $j_k$ of a basic path $q$ be denoted by $i^q$ and, respectively, $j^q$. The pair $(i^q, j^q)$ could be an activity (an element of $\mathcal{E}_{\mathrm{a}}$), in which case it is necessarily a nonbasic activity (i.e., an element of $\mathcal{E}_{\mathrm{a}} \setminus \mathcal{E}_{\mathrm{ba}}$), and we say that the graph $(\mathcal{V}_q, \mathcal{E}_q \cup \{(i^q, j^q)\})$ is a *closed simple path*; otherwise $(i^q, j^q)$ is not an activity (i.e., it is in $\mathcal{E}_{\mathrm{a}}^c$) and we say that $q$ itself is an *open simple path*. We say that $p$ is a *simple path* if it is either a closed or an open simple path. Denote by $CSP$ and $OSP$ the sets of closed and open simple paths, respectively, and by $SP$ the set of simple paths. For a path $p \in SP$, we write $\mathcal{V}_p$ and $\mathcal{E}_p$ for its vertex and edge sets, respectively. Finally, if $p$ is a simple path, let $q^p \in BP$ denote the corresponding basic path $q$, and let $i^p \in \mathcal{I}$ and $j^p \in \mathcal{J}$ denote the leaves $i^q$ and $j^q$ of $q^p$.

Note that if $p = (\mathcal{V}_p, \mathcal{E}_p)$ is a simple path and $q^p = (\mathcal{V}_q, \mathcal{E}_q)$ is its corresponding basic path, then $\mathcal{V}_q = \mathcal{V}_p$, and $\mathcal{E}_q = \mathcal{E}_p \setminus \{(i^p, j^p)\}$.

Next, we associate directions with edges of simple paths. Let $p$ be a simple path and let $q = q^p = (\mathcal{V}_q, \mathcal{E}_q)$ be the corresponding basic path. Write $\mathcal{E}_q = \{(i_0, j_0), \ldots, (i_k, j_k)\}$, where $i_0, \ldots, i_k \in \mathcal{I}$ and $j_0, \ldots, j_k \in \mathcal{J}$. The direction that will be associated with the edges in $\mathcal{E}_q$, when considered as edges of $p$, is as follows: $j_k \to i_k \to j_{k-1} \to i_{k-1} \to \cdots \to j_0 \to i_0$. In the case of an open simple path, this exhausts all edges of $p$. In the case of a closed simple path, the direction of $(i^p, j^p) = (i_0, j_k)$ is $i_0 \to j_k$. We note that an edge (corresponding to a basic activity) may have different directions when considered as an edge of different simple paths. We signify the directions along simple paths by $s(p, i, j)$, defined for $i \in \mathcal{I}$, $j \in \mathcal{J}$, $(i, j) \in \mathcal{E}_p$, $p \in SP$, as

$$(3.1) \qquad s(p, i, j) = \begin{cases} -1, & \text{if } (i, j), \text{ considered as an edge of } p, \\ & \qquad \text{is directed from } i \text{ to } j, \\ 1, & \text{if } (i, j), \text{ considered as an edge of } p, \\ & \qquad \text{is directed from } j \text{ to } i. \end{cases}$$

We will denote

$$(3.2) \qquad \mu(p) = \sum_{(i,j) \in \mathcal{E}_p} s(p, i, j) \mu_{ij}, \qquad i \in \mathcal{I}.$$

THEOREM 2. *Let Assumptions 1 and 2 hold. Then the following statements are equivalent:*

1. *The static fluid model is throughput suboptimal.*
2. *There exists a simple path $p \in SP$ such that $\mu(p) < 0$.*

Condition (2.13) is stated in terms of the variables $\{\xi_{ij}\}$. It will be convenient to work with the variables $\{\psi_{ij}\}$ in the proof below. To this end, recall that $\nu_j > 0$ for all $j$ and $\psi_{ij}^* = \xi_{ij}^* \nu_j$. Thus, the *negation* of (2.13) can be written as follows: *There exists*

$$(3.3) \qquad \psi = (\psi_{ij})_{j \in \mathcal{J}}^{i \in \mathcal{I}}, \qquad \psi_{ij} \in \mathbb{R}_+, i \in \mathcal{I}, j \in \mathcal{J},$$

*satisfying*

$$(3.4) \qquad \begin{array}{ll} \text{(a)} & \sum_{i \in \mathcal{I}} \psi_{ij} \leq \nu_j \qquad \text{for all } j \in \mathcal{J}, \\[2mm] \text{(b)} & \sum_{j \in \mathcal{J}} \psi_{ij} \leq x_i^* \qquad \text{for all } i \in \mathcal{I}, \\[2mm] \text{(c)} & \sum_{(i,j) \in \mathcal{E}} \mu_{ij} \psi_{ij} > \sum_{i \in \mathcal{I}} \lambda_i. \end{array}$$

PROOF THAT STATEMENT 2 OF THEOREM 2 IMPLIES STATEMENT 1. Assume that statement 2 holds and fix a simple path $p$ with $\mu(p) < 0$. Let $q = q^p$ be the corresponding basic path, and recall that $\psi_{ij}^* > 0$ for $(i, j) \in \mathcal{E}_q$. Denote

$$(3.5) \qquad \alpha = \min_{(i,j) \in \mathcal{E}_q} \psi_{ij}^* > 0.$$

For each $(i, j) \in \mathcal{I} \times \mathcal{J}$ we define

$$(3.6) \qquad \sigma_{ij} = -\alpha s(p, i, j) \qquad \text{if } (i, j) \in \mathcal{E}_p,$$

and $\sigma_{ij} = 0$ otherwise. Let

$$(3.7) \qquad \psi_{ij} = \psi_{ij}^* + \sigma_{ij} \qquad \text{for } (i,j) \in \mathcal{I} \times \mathcal{J}.$$

To show that statement 1 of the theorem holds, let us show that $\psi$ satisfies (3.3) and (3.4). For $(i,j) \notin \mathcal{E}_p$, $\psi_{ij} = \psi_{ij}^* \geq 0$. For $(i,j) \in \mathcal{E}_q$,

$$\psi_{ij} \geq \psi_{ij}^* - \alpha \geq 0,$$

by (3.5). In the case where $p$ is an open simple path $\mathcal{E}_p = \mathcal{E}_q$, and (3.3) follows. In the case where $p$ is a closed simple path, it is left to show that $\psi_{i^p j^p} \geq 0$. Recall that the direction associated with $(i^p, j^p)$ is $i^p \to j^p$. Thus, by (3.1) and (3.6), $\psi_{i^p j^p} \geq \psi_{i^p j^p}^* = 0$, establishing (3.3).

Next, if $p$ is a closed simple path then every vertex $v$ of it has exactly two neighbors along $p$, say, $v'$ and $v''$, and the directions of the corresponding edges are $v' \to v$ and $v \to v''$. Hence, by (3.1) and (3.6), $\sum_{j \in \mathcal{J}} \sigma_{ij} = 0$ holds for every $i \in \mathcal{I}$, and $\sum_{i \in \mathcal{I}} \sigma_{ij} = 0$ holds for every $j \in \mathcal{J}$. The term $\sigma_{i^p j^p}$, which is positive in the case when $p$ is closed, is in fact zero in the case when $p$ is open, thus yielding $\sum_{j \in \mathcal{J}} \sigma_{ij} \leq 0$ for every $i \in \mathcal{I}$ and $\sum_{i \in \mathcal{I}} \sigma_{ij} \leq 0$ for every $j \in \mathcal{J}$. Since $\psi^*$ satisfies (3.4)(a) and (b), it follows from (3.7) that so does $\psi$. Finally, by (3.7) and since (3.4)(c) holds for $\psi^*$ with equality, it suffices to prove

$$(3.8) \qquad \sum_{(i,j) \in \mathcal{E}} \mu_{ij} \sigma_{ij} > 0$$

to establish that $\psi$ satisfies (3.4)(c). By (3.6) and (3.2),

$$\sum_{(i,j) \in \mathcal{E}} \mu_{ij} \sigma_{ij} = -\alpha \sum_{(i,j) \in \mathcal{E}_p} \mu_{ij} s(p, i, j) = -\alpha \mu(p) > 0,$$

where the inequality follows from (3.5) and the assumption $\mu(p) < 0$. This establishes (3.8) and completes the proof that statement 2 implies statement 1. $\square$

In the rest of this section, we prove that statement 1 of the theorem implies statement 2. Define

$$(3.9) \qquad M(\sigma) := \sum_{(i,j) \in \mathcal{E}} \mu_{ij} \sigma_{ij},$$

for any matrix $\sigma \in \mathbb{R}^{\mathcal{E}}$. Let $S$ denote the set of $\sigma \in \mathbb{R}^{\mathcal{E}}$ satisfying the conditions

$$(3.10) \qquad \sum_{j \in \mathcal{J}} \sigma_{ij} \leq 0 \qquad \text{for all } i \in \mathcal{I}, \qquad \sum_{i \in \mathcal{I}} \sigma_{ij} \leq 0 \qquad \text{for all } j \in \mathcal{J},$$

$$(3.11) \quad \psi_{ij}^* + \sigma_{ij} \geq 0 \qquad \text{for all } (i,j) \in \mathcal{E}_a$$

and

(3.12)                    $\sigma_{ij} = 0$        for all $(i,j) \in \mathcal{E}_{\mathrm{a}}^c$.

Note that $S$ is nonempty and compact, and let

(3.13)                    $M_{\max} = \max\{M(\sigma) : \sigma \in S\}$.

The following is straightforward.

PROPOSITION 1.  *Under the assumptions of Theorem 2, the condition* $M_{\max} > 0$ *is equivalent to the throughput suboptimality of the static fluid model.*

REMARK 2.  The above result is useful in checking whether the static fluid model is throughput suboptimal, by simply solving the linear program (3.10)–(3.13) and checking if $M_{\max} > 0$.

Throughout what follows, let statement 1 hold. Proposition 1 implies $M_{\max} > 0$. Let $S_{\mathrm{opt}}$ (resp., $S_+$) denote the set of $\sigma \in S$ such that $M(\sigma) = M_{\max}$ [resp., $M(\sigma) > 0$], and note that $S_{\mathrm{opt}}$ and $S_+$ are nonempty. For $\sigma \in S_+$, consider the graph $G_\sigma = (V_\sigma, E_\sigma)$, where

(3.14)                    $E_\sigma = \{(i,j) \in \mathcal{E}_{\mathrm{a}} : \sigma_{ij} \neq 0\}$

and $V_\sigma = \{i \in \mathcal{I} : (i,j) \in E_\sigma \text{ some } j\} \cup \{j \in \mathcal{J} : (i,j) \in E_\sigma, \text{ some } i\}$ consists of all corresponding vertices. Since $M(\sigma) > 0$, we have:

(3.15)                    there exists $(i,j) \in E_\sigma$ with $\sigma_{ij} > 0$.

By (3.10) and (3.14),

(3.16)                    if $(i,j)$ is a leaf edge of $G_\sigma$ then $\sigma_{ij} < 0$,

and

(3.17)         if $(i,j) \in E_\sigma$ and $\sigma_{ij} > 0$ then there exist two edges
               $(i,j_0), (i_0,j) \in E_\sigma$ with $\sigma_{i_0,j} < 0$ and $\sigma_{i,j_0} < 0$.

DEFINITION 3.  Let $\sigma \in S_+$ be given. A subgraph $g = (\mathcal{V}_g, \mathcal{E}_g)$ of the graph $G_\sigma$ is called a *good path for* $\sigma$, if it satisfies the following conditions.

   (i)  (Connectivity) All vertices in $\mathcal{V}_g$ communicate via the edges in $\mathcal{E}_g$.
   (ii)  The degree of each vertex is at most 2.
   (iii)  The number of edges is at least 3.
   (iv)  (Alternating signs) Whenever $(i_1,j),(i_2,j) \in E_\sigma$, one has $\sigma_{i_1,j}\sigma_{i_2,j} < 0$; whenever $(i,j_1),(i,j_2) \in E_\sigma$, one has $\sigma_{i,j_1}\sigma_{i,j_2} < 0$.

(v) (Maximality) Whenever $g$ is a subgraph of some subgraph $g'$ of $G_\sigma$, and $g'$ satisfies properties (i)–(iv) above, one has $g' = g$.

It is not hard to see that observations (3.15) and (3.17) about the graph $G_\sigma$ imply that whenever $\sigma \in S_+$, there exists at least one good path for $\sigma$.

Let $\sigma \in S_+$ be given. For any edge $(i,j) \in E_\sigma$, define $s_\sigma(i,j) = -\operatorname{sign}(\sigma_{ij})$ and for any good path $g$ for $\sigma$ set

$$(3.18) \qquad \mu(\sigma, g) := \sum_{(i,j) \in \mathcal{E}_g} s_\sigma(i,j) \mu_{ij}.$$

We write $\bar{S}_C$ (resp., $\bar{S}_O$) for the set of all $\sigma \in S_{\mathrm{opt}}$ for which there exists a good path (resp., there exists no good path) $g$ for $\sigma$ that is a cycle. The letters $C$ and $O$ are mnemonics for closed and open. Note that $S_{\mathrm{opt}} = \bar{S}_C \cup \bar{S}_O$. We also set

$S_C = \{(\sigma, g) : \sigma \in \bar{S}_C \text{ and } g \text{ is a good path for } \sigma \text{ that is a cycle}\}$,

$S_O = \{(\sigma, g) : \sigma \in \bar{S}_O \text{ and } g \text{ is a good path for } \sigma \text{ that is not a cycle}\}$.

LEMMA 1. *Let $(\sigma, g) \in S_O$. Write $g = (\mathcal{V}_g, \mathcal{E}_g)$, where*

$$\mathcal{V}_g = \{v_1, \ldots, v_k\}, \qquad \mathcal{E}_g = \{(v_1, v_2), \ldots, (v_{k-1}, v_k)\},$$

*and $v_1, \ldots, v_k$ are distinct elements of $V_\sigma$. Then $\sigma_{v_1, v} < 0$ for every edge $(v_1, v) \in E_\sigma$, and similarly $\sigma_{v, v_k} < 0$ for $(v, v_k) \in E_\sigma$.*

PROOF. Argue by contradiction and assume that $\sigma_{v_1, v_2} > 0$. By (3.17), $(v_1, v_2)$ must have a neighbor $(v_0, v_1)$ with $v_0 \neq v_2$, satisfying $\sigma_{v_0, v_1} < 0$. It is easy to see that if we had $v_0 \in \mathcal{V}_g$, there would exist a good path for $\sigma$ that is a cycle, violating the assumption of the lemma that there exist no such good paths for $\sigma$. Define a new graph $g'$ by $\mathcal{V}_{g'} = \mathcal{V}_g \cup \{v_0\}$ and $\mathcal{E}_{g'} = \mathcal{E}_g \cup (v_0, v_1)$, and note that it is a good path (cf. Definition 3). This contradicts the maximality property [Definition 3(v)] and, therefore, one must have $\sigma_{v_1, v_2} < 0$. The second leaf edge $(v_{k-1}, v_k)$ is treated similarly.

To prove the second statement of the lemma, let $v_0 \neq v_2$ be such that $v_0 \in V_\sigma$, $(v_0, v_1) \in E_\sigma$. Argue by contradiction and assume that $\sigma_{v_0, v_1} > 0$. Since we already proved that $\sigma_{v_1, v_2} < 0$, we can again use the assumption that there is no good path for $\sigma$ that is a cycle to conclude that $v_0 \notin V_\sigma$. Defining $g'$ by $\mathcal{V}_{g'} = \mathcal{V}_g \cup \{v_0\}$ and $\mathcal{E}_{g'} = \mathcal{E}_g \cup (v_0, v_1)$ produces a good path that contains $g$, contradicting property (v) of Definition 3. Hence, $\sigma_{v_0, v_1} < 0$. □

LEMMA 2. *Let $(\sigma, g) \in S_C \cup S_O$. Then there exists a set $SP_g \subset SP$ of simple paths, such that*

$$(3.19) \qquad \mu(\sigma, g) = \sum_{p \in SP_g} \mu(p).$$

PROOF. Consider first the case where $(\sigma, g) \in S_C$. Write $g = (\mathcal{V}_g, \mathcal{E}_g)$ and let $\gamma^0 \in \mathbb{R}^{\mathcal{E}}$ be defined by

$$(3.20) \qquad \gamma^0_{ij} = \begin{cases} \text{sign}(\sigma_{ij}), & \text{if } (i,j) \in \mathcal{E}_g, \\ 0, & \text{otherwise.} \end{cases}$$

By (3.18) and (3.20), we have

$$(3.21) \qquad M(\gamma^0) \equiv \sum_{(i,j) \in \mathcal{E}} \gamma^0_{ij} \mu_{ij} = \sum_{(i,j) \in \mathcal{E}_g} \text{sign}(\sigma_{ij}) \mu_{ij} = -\mu(\sigma, g).$$

The following property is due to (3.20) and the fact that $g$ is a good path for $\sigma$ that is a cycle:

$$(3.22) \quad \text{for any } i \in \mathcal{I} \text{ and } j \in \mathcal{J} \text{ we have } \sum_{j \in \mathcal{J}} \gamma^0_{ij} = 0 \text{ and } \sum_{i \in \mathcal{I}} \gamma^0_{ij} = 0.$$

Define a finite sequence $\gamma^r \in \mathbb{R}^{\mathcal{E}}$ recursively as follows. Given $\gamma^r$, if there are no nonbasic activities (i.e., elements of $\mathcal{E}_a \setminus \mathcal{E}_{ba}$) in the set of edges where $\gamma^r$ is supported then terminate, and set $R = r$. Otherwise, select such a nonbasic activity, and let $p_r$ denote the (unique) closed simple path containing it as an edge. Define $\gamma^{r+1} \in \mathbb{R}^{\mathcal{E}}$ by

$$(3.23) \qquad \gamma^{r+1}_{ij} = \begin{cases} \gamma^r_{ij} + s(p_r, i, j), & \text{if } (i,j) \in \mathcal{E}_{p_r}, \\ \gamma^r_{ij}, & \text{otherwise.} \end{cases}$$

For $0 \le r < R$, the selected nonbasic activity at step $r$ is $(i^{p_r}, j^{p_r})$ (using the notation from Definition 2). By the discussion following Definition 2, the direction for this activity is $i^{p_r} \to j^{p_r}$, and thus by (3.1), we have that

$$(3.24) \qquad \gamma^{r+1}_{i,j} = \gamma^r_{i,j} - 1, \qquad \text{where } (i,j) = (i^{p_r}, j^{p_r}).$$

Given a nonbasic activity $(i,j)$, let $r$ be the first $r'$ for which $(i,j)$ is the selected nonbasic activity at step $r'$ (if such $r'$ exists). Since the transformation (3.23) modifies $\gamma$ only at basic activities and at the nonbasic activity selected at the given step, it follows that $\gamma^r_{i,j} = \gamma^0_{i,j}$. Hence, by (3.11), (3.14) and (3.20) that $\gamma^0_{i,j} = 1$. Thus, (3.24) shows that $\gamma^{r+1}_{i,j} = 0$. As a result, the support of $\gamma^{r+1}$ contains one nonbasic activity less that that of $\gamma^r$. It follows that $R < \infty$. Thus, $\gamma^R$ is well defined and supported on basic activities.

Next, since by construction, the selected simple paths are closed, it follows by the linearity of the transformation (3.23) that (3.22) holds for each $\gamma^r$, and in particular, for $\gamma^R$. It also follows from the linearity of (3.23), using (3.2), that $M(\gamma^{r+1}) = M(\gamma^r) + \mu(p_r)$ for $0 \le r < R$. Hence,

$$(3.25) \qquad M(\gamma^R) = M(\gamma^0) + \sum_{r=0}^{R-1} \mu(p_r).$$

The fact that $\gamma^R$ is supported on basic activities and that these form a tree (cf. Assumption 2), combined with the fact that $\gamma^R$ satisfies (3.22) implies that $\gamma^R = 0$. Hence, $M(\gamma^R) = 0$, and using (3.21) and (3.25), we obtain (3.19).

Next, consider $(\sigma, g) \in S_O$. Then $g = (\mathcal{V}_g, \mathcal{E}_g)$, where

$$\mathcal{V}_g = \{v_1, \dots, v_k\}, \qquad \mathcal{E}_g = \{(v_1, v_2), \dots, (v_{k-1}, v_k)\},$$

and $v_1, \dots, v_k$ are distinct elements of $V_\sigma$. First, note that either $(v_1, v_k) \in \mathcal{I} \times \mathcal{J}$ or $(v_k, v_1) \in \mathcal{I} \times \mathcal{J}$. Indeed, by Lemma 1 and properties (iii)–(iv) of Definition 3, $|\mathcal{E}_g|$ is an odd number, while having both $v_1$ and $v_k$ belong to either $\mathcal{I}$ or $\mathcal{J}$ would result with an even number for $|\mathcal{E}_g|$.

Also, we claim that $(v_1, v_k) \in \mathcal{E}_a^c$. Argue by contradiction and assume that $\mu_{v_1, v_k} > 0$. If we had $\sigma_{v_1, v_k} > 0$, then by Lemma 1, $\sigma_{v_1, v_2} < 0$ and $\sigma_{v_{k-1}, v_k} < 0$, and there would exist a good path, which is a cycle. This is prohibited since $(\sigma, g) \in S_O$. Hence, $\sigma_{v_1, v_k} \le 0$. Set $\delta := \min\{|\sigma_{v_1, v_2}|, |\sigma_{v_{k-1}, v_k}|\}$ and define a new matrix $\sigma' \in \mathbb{R}^{\mathcal{E}}$ by assigning $\sigma'_{v_1, v_k} = \sigma_{v_1, v_k} + \delta$ and $\sigma'_{ij} = \sigma_{ij}$ for $(i, j) \in \mathcal{E} \setminus \{(v_1, v_k)\}$. By the definition $\sigma'$ satisfies (3.10)–(3.12) (see also Lemma 1), which implies $M(\sigma') = M(\sigma) + \delta \mu_{v_1, v_k} > M(\sigma)$. This contradicts the assumption $\sigma \in S_{\text{opt}}$. Therefore, $\mu_{v_1, v_k} = 0$ meaning $(v_1, v_k) \in \mathcal{E}_a^c$.

The rest of the argument is similar to the treatment of the case where $(\sigma, g) \in S_C$, with some modifications as follows. Instead of (3.20), consider $\gamma^0 \in \mathbb{R}^{\mathcal{E}}$ defined as

$$(3.26) \qquad \gamma_{ij}^0 = \begin{cases} \text{sign}(\sigma_{ij}), & \text{if } (i, j) \in \mathcal{E}_g, \\ 1, & \text{if } (i, j) = (v_1, v_k), \\ 0, & \text{otherwise.} \end{cases}$$

Since $(v_1, v_k)$ is not an activity, $\mu_{v_1, v_k} = 0$, and thus (3.21) is still valid. Also, it follows from Lemma 1 that $\gamma_{v_1, v_2}^0 = \gamma_{v_{k-1}, v_k}^0 = -1$, and as a result, (3.22) is valid. We can now repeat the construction of $\{\gamma^r\}$, $0 \le r < R$, and the inductive argument that leads to (3.25). The matrix $\gamma^R$, in this case, is supported on basic activities plus the edge $(v_1, v_k)$. Denoting by $p$ the open simple path whose leaves are $v_1$ and $v_k$, we apply one last time a transformation of the form (3.23) as follows:

$$\gamma_{ij}^{R+1} = \begin{cases} \gamma_{ij}^R + s(p, i, j), & \text{if } (i, j) \in \mathcal{E}_p, \\ \gamma_{ij}^R - 1, & \text{if } (i, j) = (i^p, j^p) \equiv (v_1, v_k), \\ \gamma_{ij}^R, & \text{otherwise.} \end{cases}$$

As a result,

$$M(\gamma^{R+1}) = M(\gamma^0) + \sum_{r=0}^{R-1} \mu(p_r) + \mu(p).$$

Arguing as before, we obtain that $\gamma^{R+1}$ is supported on basic activities, satisfies (3.22) thus vanishes. Hence, $M(\gamma^{R+1}) = 0$, and (3.19) follows as before.  □

LEMMA 3.    *Let $(\sigma, g) \in S_C \cup S_O$. Then the following statements are true.*

   (i)  $\mu(\sigma, g) \leq 0$.
   (ii) *If $\mu(\sigma, g) = 0$ then there exists $\sigma' \in S_{\mathrm{opt}}$ satisfying $E_{\sigma'} \subsetneq E_\sigma$.*

PROOF.    We first prove part (i). Consider first the case where $(\sigma, g) \in S_C$. Arguing by contradiction, assume $\mu(\sigma, g) > 0$. Let

$$(3.27) \qquad \alpha = \min_{(i,j) \in \mathcal{E}_g} |\sigma_{ij}| > 0,$$

and for each $(i, j) \in E_\sigma$ define

$$(3.28) \qquad \sigma'_{ij} = \sigma_{ij} + s_\sigma(i,j)\alpha \quad \text{if } (i,j) \in \mathcal{E}_g, \quad \text{and} \quad \sigma'_{ij} = \sigma_{ij} \quad \text{otherwise.}$$

We show that $\sigma'$ satisfies conditions (3.10)–(3.12) and, therefore, that $\sigma' \in S$. To this end, note that the sums in (3.10) remain unchanged under the transformation from $\sigma$ to $\sigma'$, due to the fact that $g$ is a cycle and using the alternating signs property [Definition 3(iv)]. Thus, (3.10) is satisfied by $\sigma'$. The relation (3.11) follows from (3.27) and (3.28), since $\sigma'_{ij} > \sigma_{ij}$ for $(i, j) \in \mathcal{E}_g$ with $\sigma_{ij} < 0$ and $\sigma'_{ij} \geq 0$ for $(i, j) \in \mathcal{E}_g$ with $\sigma_{ij} > 0$. The relation (3.12) holds trivially. This shows $\sigma' \in S$. We have

$$(3.29) \qquad \begin{aligned} M(\sigma') &= \sum_{(i,j) \in \mathcal{E}} \sigma'_{ij}\mu_{ij} = \sum_{(i,j) \in \mathcal{E}} \sigma_{ij}\mu_{ij} + \alpha \sum_{(i,j) \in \mathcal{E}_g} s_\sigma(i,j)\mu_{ij} \\ &= M(\sigma) + \alpha\mu(\sigma, g). \end{aligned}$$

Since $\mu(\sigma, g) > 0$ by assumption, we have $M(\sigma') > M(\sigma)$, which contradicts the assumption $\sigma \in S_{\mathrm{opt}}$. Hence, (i) holds.

Consider now the case where $(\sigma, g) \in S_O$. Argue by contradiction and assume that $\mu(\sigma, g) > 0$. Define $\sigma'$ as in (3.27)–(3.28). Once again, we claim that the constraints (3.10)–(3.12) are satisfied for $\sigma'$. The argument following (3.28) applies, and it remains only to check (3.10) for the vertices $v_1$ and $v_k$. The validity of (3.10) in this case follows from (3.27), (3.28) and Lemma 1, since we have $\sigma_{v_1,v_2} < 0$, $\sigma_{v_{k-1},v_k} < 0$ and $\sigma'_{ij} \leq 0$ holds for all $(i,j) \in \mathcal{E}_g$ with $\sigma_{ij} < 0$. This shows that $\sigma' \in S$. The rest of the argument is as in the previous case.

Next we prove part (ii). The desired $\sigma'$ is, in fact, the one constructed in the proof of part (i). Indeed, we have proved that $\sigma' \in S$. Moreover, since $\mu(\sigma, g) = 0$ and $\sigma \in S_{\mathrm{opt}}$, we have by (3.29) that $\sigma' \in S_{\mathrm{opt}}$. By (3.27) and (3.28),

$$\sigma'_{i'j'} = 0 \qquad \text{for all } (i', j') \in \underset{(i,j) \in \mathcal{E}_g}{\arg\min} |\sigma_{ij}|,$$

and, therefore, statement (ii) holds. □

PROOF THAT STATEMENT 1 OF THEOREM 2 IMPLIES STATEMENT 2. Let $\sigma \in S_{\mathrm{opt}}$. Let $g$ be such that $(\sigma, g) \in S_C \cup S_O$. Set $(\sigma^0, g^0) = (\sigma, g)$ and define a finite sequence $(\sigma^r, g^r) \in S_C \cup S_O$ recursively as follows. If $\mu(\sigma^r, g^r) < 0$ then set $R = r$ and terminate. Otherwise, by Lemma 3(i), $\mu(\sigma^r, g^r) = 0$. Let $\sigma^{r+1}$ denote the matrix $\sigma'$ from Lemma 3(ii) corresponding to $(\sigma^r, g^r)$. Since $\sigma^{r+1} \in S_{\mathrm{opt}} = \bar{S}_C \cup \bar{S}_O$, it follows from the definition of $S_O$ and $S_C$ that there exists $g$ such that $(\sigma^{r+1}, g) \in S_C \cup S_O$. Let $g^{r+1}$ be such a good path.

The finiteness of $R$ follows from Lemma 3(ii) and the finiteness of the set $\mathcal{E}_\sigma$.

By construction, $(\sigma^R, g^R) \in S_C \cup S_O$ and $\mu(\sigma^R, g^R) < 0$. Lemma 2 thus implies that there exists a simple path $p$ such that $\mu(p) < 0$. This concludes the proof of the theorem. □

We end this section by revisiting the three examples from Section 2.3. We can now use Theorem 2 to determine throughput suboptimality for each example.

EXAMPLE 1. The simple path $p$ corresponding to Example 1 of Section 2.3 (see Figure 2, left) satisfies $\mu(p) = -4 < 0$. Hence Theorem 1 applies. Moreover, $p$ is a *closed* simple path, and one checks that [3], Theorem 2.3, is valid, also.

EXAMPLE 2. In the case of Example 2 of Section 2.3, the simple path $p$ is open (see Figure 2, right). We have $\mu(p) = -3 < 0$. Theorem 1 applies. Since $p$ is open, [3], Theorem 2.3, does not apply.

EXAMPLE 3. To see that the static fluid model of Example 3 of Section 2.3 is throughput optimal, we use Proposition 1 and solve the linear optimization problem (3.10)–(3.13), to find that $M_{\max} = 0$. Hence, throughput optimality holds.

## 4. Dynamic fluid model.

As a tool for analyzing the probabilistic model, we consider a model with deterministic arrival and service rates, and deterministic perturbations. The model is a deterministic analogue of a model that will be studied in Section 5, which corresponds to our original probabilistic model at fluid scale [see equations (5.3)–(5.9)]. The deterministic model that we now introduce will be referred to as a dynamic fluid model. It consists of cadlag functions $X_i$, $Y_i$, $Z_j$, $\Psi_{ij}$, $W_i$, $i \in \mathcal{I}$, $j \in \mathcal{J}$, satisfying

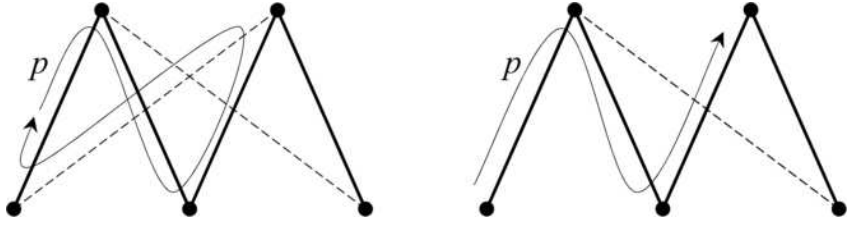

FIG. 2.  *Simple paths for Examples* 1 *and* 2: *On the left p is a closed simple path, while on the right p is open. For Example* 1, $\mu_{21} > 0$ *and* $(2,1)$ *is a nonbasic activity. For Example* 2, $\mu_{21} = 0$ *and* $(2,1)$ *is not an activity.*

the equations below, for all $t \geq 0$.

$$(4.1) \qquad X_i(t) = x_i^* + W_i(t) + \lambda_i t - \sum_{j \in \mathcal{J}} \mu_{ij} \int_0^t \Psi_{ij}(s)\, ds, \qquad i \in \mathcal{I},$$

$$(4.2) \quad Y_i(t) + \sum_{j \in \mathcal{J}} \Psi_{ij}(t) = X_i(t), \qquad i \in \mathcal{I},$$

$$(4.3) \quad Z_j(t) + \sum_{i \in \mathcal{I}} \Psi_{ij}(t) = \nu_j + \theta_j, \qquad j \in \mathcal{J},$$

$$(4.4) \qquad \Psi_{ij}(t) = 0, \qquad (i,j) \in \mathcal{E}_{\mathrm{a}}^c,$$

$$(4.5) \qquad Y_i(t) \geq 0, \qquad Z_j(t) \geq 0, \qquad \Psi_{ij}(t) \geq 0, \qquad i \in \mathcal{I},\ j \in \mathcal{J}.$$

Above, the constants $\{x_i^*, \lambda_i, \mu_{ij}, \nu_j; i \in \mathcal{I}, j \in \mathcal{J}\}$ are as in Section 2, and $\{\theta_j; j \in \mathcal{J}\}$ are additional real constants. We refer to $(W, \theta)$ as *data* for the model. The term $W_i$ is analogous to the process (5.8) of the probabilistic model, representing perturbations in the arrival and service processes that are eventually shown to converge to zero. The term $\theta$ corresponds to a discretization error due to the fact that the quantities $N_j^n$ take integer values.

DEFINITION 4.  Given $\sigma > \varepsilon > 0$, we will say that the data $(W, \theta)$ for the dynamic fluid model (4.1)–(4.5) is an $(\varepsilon, \sigma)$-*perturbation* if

$$(4.6) \qquad \|W(t)\| \leq \varepsilon \qquad \text{for all } 0 \leq t < \sigma \quad \text{and} \quad \|\theta\| \leq \varepsilon.$$

Let Assumptions 1 and 2 hold and, in addition, assume that the static fluid model is throughput suboptimal. Let $\varepsilon > 0$ and $\sigma > 0$ be given and assume that the data $(W, \theta)$ is an $(\varepsilon, \sigma)$-perturbation. Below we will construct functions $\{X_i(t), i \in \mathcal{I}\}$, $\{Y_i(t), i \in \mathcal{I}\}$, $\{Z_j(t), j \in \mathcal{J}\}$ and $\{\Psi_{ij}(t), (i,j) \in \mathcal{I} \times \mathcal{J}\}$ that satisfy (4.1)–(4.5), such that $\int_0^\sigma 1_{\{e \cdot Y(s) > 0\}}\, ds$ and $\|X(t) - x^*\|_\sigma^*$ are $o(1)$ as $\varepsilon$ becomes small. The precise statement will be formulated in Theorem 3.

Since the static fluid model is throughput suboptimal, we have from Theorem 2 that there exists a simple path $p$ with $\mu(p) < 0$. Fix such a path $p$.

Set

$$\mathcal{E}_p^+ = \{(i,j) \in \mathcal{E}_p : s(p,i,j) = 1\}, \qquad \mathcal{E}_p^- = \{(i,j) \in \mathcal{E}_p : s(p,i,j) = -1\},$$

and note that

$$(4.7) \qquad \mu(p) = \sum_{(i,j) \in \mathcal{E}_p} s(p,i,j)\mu_{ij} = \Sigma_p^+ - \Sigma_p^- < 0,$$

where

$$\Sigma_p^+ := \sum_{(i,j) \in \mathcal{E}_p^+} \mu_{ij}, \qquad \Sigma_p^- := \sum_{(i,j) \in \mathcal{E}_p^-} \mu_{ij}.$$

In addition, set

$$(4.8) \qquad \Sigma_p^0 := \sum_{(i,j) \in \mathcal{E}_p^c \cap \mathcal{E}_{\mathrm{ba}}} \mu_{ij}.$$

Define the constant

$$(4.9) \qquad \alpha = \tfrac{1}{2}(1 + \Sigma_p^+ / \Sigma_p^-).$$

The following inequalities follow from (4.7) and (4.9)

$$(4.10) \qquad \tfrac{1}{2} < \alpha < 1, \qquad \Sigma_p^+ - \alpha \Sigma_p^- = \tfrac{1}{2}(\Sigma_p^+ - \Sigma_p^-) < 0.$$

We will need the following result from [2], Proposition A.2 (that follows from the tree structure of $\mathcal{E}_{\mathrm{ba}}$): the system of equations

$$(4.11) \qquad \begin{cases} \displaystyle\sum_{j \in \mathcal{J}} \phi_{ij} = a_i, & i \in \mathcal{I}, \\ \displaystyle\sum_{i \in \mathcal{I}} \phi_{ij} = b_j, & j \in \mathcal{J}, \\ \phi_{ij} = 0, & (i,j) \in \mathcal{E}_{\mathrm{ba}}^c, \end{cases}$$

in the unknown $\phi$ has a unique solution, whenever $a$ and $b$ satisfy $\sum_i a_i = \sum_j b_j$. With

$$(4.12) \qquad D_G = \left\{ (a,b) \in \mathbb{R}^I \times \mathbb{R}^J : \sum_{i \in \mathcal{I}} a_i = \sum_{j \in \mathcal{J}} b_j \right\},$$

denote by $G : D_G \to \mathbb{R}^{I \times J}$ the solution map, namely

$$(4.13) \qquad \phi_{ij} = G_{ij}(a,b), \qquad (i,j) \in \mathcal{E},$$

and note that this map is linear. Let

$$(4.14) \qquad C_G := \sup \left\{ \max_{ij} |G_{ij}(a,b)| : (a,b) \in D_G, \|a\| \vee \|b\| \le 1 \right\}.$$

Recall that by the definition of basic activities, $\psi^*$ has the property that $\psi^*_{ij} > 0$ for $(i,j) \in \mathcal{E}_{ba}$. Set

$$(4.15) \qquad \delta_1 = \tfrac{1}{2} \min_{(i,j) \in \mathcal{E}_{ba}} \psi^*_{ij} > 0, \qquad a_0 = (2C_G)^{-1}\delta_1,$$

$$(4.16) \qquad \delta_2 = \begin{cases} \delta_1 \min\left\{ \dfrac{\alpha\Sigma^-_p - \Sigma^+_p}{2\Sigma^0_p}, (1-\alpha) \right\}, & \text{if } \Sigma^0_p > 0, \\ 0, & \text{if } \Sigma^0_p = 0. \end{cases}$$

By (4.10), we have $0 \le \delta_2 < \delta_1$.

Recall that there exists a simple path $p$ with $\mu(p) < 0$ and let such $p$ be fixed. Let $\beta = \beta(\varepsilon) \in \mathbb{R}^{\mathcal{E}}$ be a constant matrix satisfying

$$(4.17) \qquad |\beta| := \max_{(i,j) \in \mathcal{E}} |\beta_{ij}| \le \varepsilon^2.$$

Define the constant matrix $\tilde{\psi} \in \mathbb{R}^{\mathcal{E}}$ as

$$(4.18) \qquad \tilde{\psi}_{ij} = \begin{cases} \psi^*_{ij} + \alpha\delta_1 + \beta_{ij}, & (i,j) \in \mathcal{E}^-_p, \\ \psi^*_{ij} - \delta_1 + \beta_{ij}, & (i,j) \in \mathcal{E}^+_p, \\ \psi^*_{ij} - \delta_2 + \beta_{ij}, & (i,j) \in \mathcal{E}^c_p \cap \mathcal{E}_{ba}, \\ 0, & (i,j) \in \mathcal{E}^c_p \cap \mathcal{E}^c_{ba}. \end{cases}$$

The term $\beta$ will be required in Section 5 where the results are translated to the probabilistic model, in which the process $\Psi^n$ takes integer values.

Throughout, we use $X_e$ for $e \cdot X$, and use a similar convention for $x^*_e, \nu_e$ and $\theta_e$.

We now explain the main idea of the construction. Recall that our goal is to come up with functions $X, Y, Z, \Psi$ satisfying the dynamic fluid model equations, where in addition, (i) the function $Y$ takes the value zero "most of the time," and (ii) the function $X$ is kept close to the quantity $x^*$ of the static fluid model. The construction will be based on alternating between two different behaviors. On intervals that will be denoted by $[\eta_{k-1}, \zeta_k)$, equation (4.20) will be in force. These intervals will be short, and the function $Y$ will not necessary vanish on them, but they will achieve a quick reduction of the value of $X_e$, so that at the end of such an interval, $X_e < 0$. As a result, at the end of the interval, one will be able to find $\Psi$ with $Y = 0$ (because $X_e < 0$ corresponds to the number of customers in the system being smaller than the total number of servers). The quick reduction of $X_e$ will be achieved by using an allocation matrix under which the total processing rate is greater than the total arrival rate (2.11), which exists, by throughput suboptimality. On intervals $[\zeta_k, \eta_k)$, equation (4.22) will be valid. These intervals will be relatively long. The deviation of $\Psi$ from the value $\psi^*$ of the static fluid model will be small, and as we shall prove this will achieve the closeness of $X$ to $x^*$. It will also be shown that the function $Y$, that starts at zero, is kept zero on these intervals. The combination of the shortness of the intervals $[\eta_{k-1}, \zeta_k)$ and the vanishing of $Y$ on the intervals $[\zeta_k, \eta_k)$ will finally yield the desired result regarding $Y$ being zero most of the time (4.27). The construction of the functions is carried out in Lemma 4. In the proof of Theorem 3, the properties alluded to above are established.

Fix $i_0 \in \mathcal{I}$ and $j_0 \in \mathcal{J}$. Set $\eta_0 = 0$. Let $k \geq 1$ and consider the system of equations (4.19)–(4.23):

$$(4.19) \qquad \tau = \tilde{\tau} \wedge \sigma, \qquad \tilde{\tau} = \inf\{t \geq 0, \|X(t) - x^*\| \geq \varepsilon^{1/2}\},$$

$$(4.20) \quad
\begin{cases}
X_i(t) = X_i(\eta_{k-1}) + W_i(t) - W_i(\eta_{k-1}) + \lambda_i(t - \eta_{k-1}) \\
\qquad\qquad - \sum_{j \in \mathcal{J}} \mu_{ij} \int_{\eta_{k-1}}^t \Psi_{ij}(s)\, ds, \qquad i \in \mathcal{I},\ t \in [\eta_{k-1}, \zeta_k), \\
\Psi(t) = \tilde{\psi}, \qquad\qquad\qquad\qquad\qquad t \in [\eta_{k-1}, \zeta_k), \\
Y_i(t) = X_i(t) - \sum_{j \in \mathcal{J}} \Psi_{ij}(t), \qquad i \in \mathcal{I},\ t \in [\eta_{k-1}, \zeta_k), \\
Z_j(t) = \nu_j + \theta_j - \sum_{i \in \mathcal{I}} \Psi_{ij}(t), \qquad j \in \mathcal{J},\ t \in [\eta_{k-1}, \zeta_k),
\end{cases}$$

where

$$(4.21) \qquad \zeta_k = \inf\{t \geq \eta_{k-1} : X_e(t) - X_e(\eta_{k-1}) \leq -7\varepsilon\} \wedge \tau,$$

and

$$(4.22) \quad
\begin{cases}
X_i(t) = X_i(\zeta_k) + W_i(t) - W_i(\zeta_k) + \lambda_i(t - \zeta_k) \\
\qquad\qquad - \sum_{j \in \mathcal{J}} \mu_{ij} \int_{\zeta_k}^t \Psi_{ij}(s)\, ds, \qquad i \in \mathcal{I}, t \in [\zeta_k, \eta_k), \\
Y(t) = (X_e(t) - \nu_e - \theta_e)^+ e_{i_0}, \\
\qquad Z(t) = (X_e(t) - \nu_e - \theta_e)^- e_{j_0}, \qquad t \in [\zeta_k, \eta_k), \\
\Psi(t) = G(X(t) - Y(t), \nu + \theta - Z(t)), \qquad t \in [\zeta_k, \eta_k),
\end{cases}$$

where

$$(4.23) \qquad \eta_k = \inf\{t \geq \zeta_k : \|X(t) - X(\zeta_k)\| \geq 3\varepsilon\} \wedge \tau.$$

LEMMA 4. *Assume that the data $(W, \theta)$ for the dynamic fluid model is an $(\varepsilon, \sigma)$-perturbation. Then for all $\varepsilon > 0$ sufficiently small, equations (4.19)–(4.23) uniquely define a number $\tau \in (0, \sigma]$, a finite sequence $\eta_k, \zeta_k$ and functions $X, Y, Z$ and $\Psi$ on $t \in [0, \tau)$. Moreover, these functions satisfy (4.1)–(4.5) on $[0, \tau)$.*

PROOF. In what follows, the functions $X, Y, Z$ and $\Psi$ are constructed recursively. The fact that these functions and the number $\tau$, satisfying (4.19), are uniquely determined by the data will follow by construction. The positivity of $\tau$ will follow from the bound $\zeta_1 \geq \varepsilon^2$, proved in the paragraph that follows.

Note next that if $k \geq 1$ and $\eta_{k-1}$ and $X(\eta_{k-1})$ are given, then the first two equations in (4.20) define $X$ and $\Psi$ uniquely on $[\eta_{k-1}, \zeta_k)$. The last two lines of (4.20) define $Y$ and $Z$ on the same interval. Moreover, it is easy to see that equations (4.1)–(4.4) are satisfied on $[0, \zeta_k)$, provided that they are satisfied on $[0, \eta_{k-1})$. The validity of the constraints (4.5) is argued later. Next, let $\zeta_k$ and $X(\zeta_k)$ be given. Substituting into the first equation in (4.22) the values for $Y, Z$ and $\Psi$ from the last three equations of (4.22) results with an equation of the form

$$X_i(t) = X_i(\zeta_k) + W_i(t) - W_i(\zeta_k) + \int_{\zeta_k}^t F_i(X(s)) \, ds, \qquad i \in \mathcal{I},$$

where $F_i \colon \mathbb{R}^I \to \mathbb{R}$ are globally Lipschitz. This uniquely defines $X$, and in turn, $Y, Z$ and $\Psi$ on $[\zeta_k, \eta_k)$. Thus, $X, Y, Z$ and $\Psi$ are uniquely defined on $[0, \eta_k)$, provided that they are on $[0, \zeta_k)$. Let

$$(4.24) \qquad K = \inf\{k \geq 0 \colon \eta_k = \tau \text{ or } \zeta_{k+1} = \tau\}.$$

We show that $K < \infty$. To this end, observe that for $\varepsilon > 0$ sufficiently small, we have $\zeta_1 \geq \varepsilon^2 > 0$. Indeed, if $\zeta_1 = \sigma$, $K = 0$. Otherwise, (4.6), (4.3) and (4.5) imply the inequality $\Psi_{ij}(t) \leq \nu_j + \varepsilon$, and as a result, from (4.20) and (4.6), it follows that for any $0 \leq t \leq \zeta_1$,

$$(4.25) \qquad |X_e(t) - X_e(0)| \leq \|X(t) - X(0)\| \leq 2\varepsilon + (c_1 + c_2\varepsilon)t,$$

where $c_1 = \sum_i \lambda_i + \sum_{ij} \mu_{ij}\nu_j$ and $c_2 = \sum_{ij} \mu_{ij}$, and the first inequality above is due to the relation $|a_e| = |a_1 + \cdots + a_d| \leq \|a\|$, $a \in \mathbb{R}^d$. By right-continuity of $X$ and (4.21), we have that $|X_e(\zeta_1) - X_e(0)| \geq 7\varepsilon$. Thus, $\zeta_1 \geq 5\varepsilon(c_1 + c_2\varepsilon)^{-1} \geq \varepsilon^2$ for $\varepsilon$ sufficiently small. For $k \leq K$, denote

$$(4.26) \qquad I_k^1 = [\eta_{k-1}, \zeta_k), \qquad I_k^2 = [\zeta_k, \eta_k), \qquad I_k = I_k^1 \cup I_k^2.$$

An argument similar to the one for $\zeta_1 \geq \varepsilon^2$ shows that each of the intervals $I_k^1$ and $I_k^2$ has a length of at least $\varepsilon^2$. Hence, $K$ is finite. The inductive argument given above thus shows that the functions $X, Y, Z$ and $\Psi$ are uniquely defined on $[0, \tau)$ and satisfy (4.1)–(4.4). We now show that the relations (4.5) hold. The two interval types are treated separately.

*Intervals $I_k^1$*: Let $k$ and $t \in I_k^1$ be fixed. By (4.20), $\Psi_{ij}(t)$ is given by $\widetilde{\psi}$, defined in (4.18). The nonnegativity of each $\widetilde{\psi}_{ij}$ for all sufficiently small $\varepsilon$ follows from (4.15) and $0 \leq \delta_2 < \delta_1$ [cf. (4.16)]. To show that $Z_j$ are nonnegative, note that by (4.20) and (4.18),

$$Z_j(t) = \nu_j + \theta_j - \sum_{i \in \mathcal{I}} \Psi_{ij}(t) = \nu_j + \theta_j - \sum_{i \in \mathcal{I}} \widetilde{\psi}_{ij}$$

$$\geq \theta_j - C|\beta| + \delta_1|\{i : (i,j) \in \mathcal{E}_p^+\}|$$
$$- \alpha\delta_1|\{i : (i,j) \in \mathcal{E}_p^-\}| + \delta_2|\{i : (i,j) \in \mathcal{E}_p^c \cap \mathcal{E}_{\mathrm{ba}}\}|,$$

where $C = |\mathcal{E}_{\mathrm{a}}|$, and we have used (2.10) in the last line. By Definition 2 and (3.1), for any vertex $j \in \mathcal{J}$ in $\mathcal{V}_p$, there exists exactly one edge, $(i_1, j) \in \mathcal{E}_p$, with $s(p, i_1, j) = 1$, and there exists at most one edge $(i_2, j) \in \mathcal{E}_p$ with $s(p, i_2, j) = -1$ (in the case where $p$ is an open simple path and $j$ is a leaf, such $i_2$ does not exist). Thus, the positivity of $\delta_1$, nonnegativity of $\delta_2$ and the bounds on $\beta$ and $\theta$ show that $Z_j(t) \geq (1 - \alpha)\delta_1 - (C + 1)\varepsilon \geq 0$ for all $\varepsilon$ sufficiently small and $j \in \mathcal{J}$.

Given $i \in \mathcal{I}$, by similar considerations,

$$\sum_{j \in \mathcal{J}} \Psi_{ij}(t) = \sum_{j \in \mathcal{J}} \widetilde{\psi}_{ij}$$
$$\leq x_i^* + C|\beta| - \delta_1|\{j : (i,j) \in \mathcal{E}_p^+\}|$$
$$+ \alpha\delta_1|\{j : (i,j) \in \mathcal{E}_p^-\}| - \delta_2|\{j : (i,j) \in \mathcal{E}_p^c \cap \mathcal{E}_{\mathrm{ba}}\}|$$
$$\leq x_i^* + C\varepsilon - (1 - \alpha)\delta_1.$$

Since $t < \tau$, we have by (4.19) that $\|X(t) - x^*\| < \varepsilon^{1/2}$. By (4.10) and (4.15), $(1 - \alpha)\delta_1 > 0$. We conclude that $\sum_{j \in \mathcal{J}} \Psi_{ij}(t) \leq X_i(t)$, and in turn by (4.20), $Y_i(t) \geq 0$, provided that $\varepsilon$ is sufficiently small.

*Intervals $I_k^2$:* Fix $k$ and $t \in I_k^2$. The nonnegativity of $Y_i(t)$ and $Z_j(t)$ is immediate from (4.22). Also, (4.22) and (4.11) imply $\Psi_{ij}(t) = 0$ for $(i,j) \in \mathcal{E}_{\mathrm{ba}}^c$. It remains to show that $\Psi_{ij}(t) \geq 0$ for $(i,j) \in \mathcal{E}_{\mathrm{ba}}$. By (4.19), $\|X(t) - x^*\| \leq a_0$ for $t < \tau$ and $\varepsilon$ sufficiently small. By uniqueness of solutions to (4.11) and by (2.10), $\psi^* = G(x^*, \nu)$. Also (2.10), (4.6), (4.19) and (4.22) imply that $\|Y(t)\| \vee \|Z(t)\| \leq \|X(t) - x^*\| + \varepsilon$ for $t < \tau$. Hence, by linearity of the map $G$ on the domain $D_G$ (4.12), by (4.22) and (4.13)–(4.15), we have

$$\Psi_{ij}(t) = G_{ij}(x^*, \nu) + G_{ij}(X(t) - x^* - Y(t), \theta - Z(t))$$
$$\geq \psi_{ij}^* - C_G(\|X(t) - x^* - Y(t)\| \vee \|\theta - Z(t)\|)$$
$$\geq \psi_{ij}^* - 2C_G(\|X(t) - x^*\| + \varepsilon)$$
$$\geq \psi_{ij}^* - 2C_G a_0 - 2C_G\varepsilon$$
$$= \psi_{ij}^* - \delta_1 - 2C_G\varepsilon$$
$$\geq 0,$$

where the last inequality holds by the definition of $\delta_1$ (4.15) and provided that $\varepsilon$ is sufficiently small. This concludes the proof of Lemma 4. $\square$

THEOREM 3. *Let Assumptions 1 and 2 hold. Assume that the static fluid model is throughput suboptimal. Then there exist functions $\gamma_1$ and $\gamma_2$ from* $(0, \infty)$ *to itself, satisfying* $\lim_{\varepsilon \to 0} \gamma_1(\varepsilon) = 0$ *and* $\lim_{\varepsilon \to 0} \gamma_2(\varepsilon) = \infty$, *such that the following statement holds. If the data* $(W, \theta)$ *for the dynamic fluid model is an* $(\varepsilon, \sigma)$-*perturbation, then the functions* $X$, $Y$, $Z$ *and* $\Psi$ *that are uniquely defined by* (4.19)–(4.23), *satisfy*

$$(4.27) \qquad \int_0^{\sigma \wedge \gamma_2(\varepsilon)} 1_{\{e \cdot Y(s) > 0\}} \, ds \le \gamma_1(\varepsilon),$$

$$(4.28) \qquad \|X(t) - x^*\| \le \gamma_1(\varepsilon) \qquad \text{for all } 0 \le t \le \sigma \wedge \gamma_2(\varepsilon).$$

The proof of the following lemma appears at the end of the section.

LEMMA 5. *Recall the definitions of* $K$ (4.24) *and intervals* $I_k^1$ *and* $I_k^2$ *from* (4.26). *There exist constants* $m_1$, $m_2$, $m_3 \in (0, \infty)$, *not depending on* $\varepsilon$, $\sigma$ *and* $k$, *such that for any* $k \le K$:

1. $|I_k^1| \le m_1 \varepsilon$;
2. $\|X - x^*\|_{\eta_k}^* \le k m_2 \varepsilon$;
3. $|I_k^2| \ge m_3 / k$.

PROOF OF THEOREM 3. We begin by showing that

$$(4.29) \qquad Y(t) = 0 \qquad \text{for all } t \in [\zeta_k, \eta_k), k < K.$$

By (4.22), it suffices to show that

$$(4.30) \qquad X_e(t) - \nu_e - \theta_e \le 0 \qquad \text{for all } t \in [\zeta_k, \eta_k), k < K.$$

Indeed, from (4.21), (4.6) and using $\nu_e = x_e^*$ [by (2.10)], we have

$$X_e(\zeta_1) - \nu_e - \theta_e \le X_e(\zeta_1) - X_e(0) + (X_e(0) - x_e^*) - \theta_e \le -7\varepsilon + \varepsilon - \theta_e \le -5\varepsilon.$$

Then by (4.23), taking into account the possibility of jumps of at most $2\varepsilon$ for $W_e$, we have for $t \in [\zeta_1, \eta_1)$:

$$X_e(t) - \nu_e - \theta_e \le X_e(\zeta_1) - \nu_e - \theta_e + \|X(t) - X(\zeta_1)\| \le -5\varepsilon + 5\varepsilon \le 0.$$

A proof by induction that repeats the above argument, using (4.21) and (4.23) shows that (4.30), and in turn (4.29), holds for $k \ge 1$.

Next, let us show that

$$(4.31) \qquad \tau \ge \sigma \wedge \gamma_2(\varepsilon),$$

where $\gamma_2(\varepsilon) := \frac{m_3}{4} |\log \varepsilon|$, for sufficiently small $\varepsilon$. Consider the number $k_0 = k_0(\varepsilon) := [(2m_2 \varepsilon^{1/2})^{-1}] \wedge K$ [where $K$ is as in (4.24)]. If $k_0 = K$ and $\tau = \eta_K$ then from (4.19), (4.24) and using Lemma 5(2), we have $\tau = \sigma$, since then $\|X - x^*\|_{\eta_K} < \varepsilon^{1/2}$. Otherwise, if $k_0 = K$ and $\tau = \zeta_{K+1}$, or $k_0 = [(2m_2\varepsilon^{1/2})^{-1}]$, one uses (4.19), (4.24) and Lemma 5(2), (3) to obtain

$$\tau \geq \eta_{k_0} \geq \sum_{l=1}^{k_0(\varepsilon)} \frac{m_3}{l} \geq \frac{m_3}{4}|\log \varepsilon|.$$

Hence, (4.31) follows.

Let $K_0 = K_0(\varepsilon) = \max\{k : \zeta_k \leq \sigma \wedge \gamma_2(\varepsilon)\}$. By Lemma 5(3), for sufficiently small $\varepsilon$,

$$\sum_{k=1}^{K_0-1} k^{-1} \leq m_3^{-1}\gamma_2(\varepsilon) = -\frac{1}{4}\log \varepsilon.$$

This implies that $\frac{1}{2}\log K_0 \leq -\frac{1}{4}\log \varepsilon$, hence, $K_0 \leq \varepsilon^{-1/2}$, provided that $\varepsilon$ is small.

Now, using (4.29), Lemma 5(1) and the estimate on $K_0$, we have

$$\int_0^{\sigma \wedge \gamma_2(\varepsilon)} 1_{\{e \cdot Y(s) > 0\}} \, ds \leq (K_0(\varepsilon) + 1)m_1\varepsilon \leq 2m_1\varepsilon^{1/2} \leq \gamma_1(\varepsilon),$$

where $\gamma_1(\varepsilon) := 2(m_1 \vee 1)\varepsilon^{1/2}$, establishing (4.27). As a result of (4.19) and (4.31), we obtain that

$$\|X(t) - x^*\| < \varepsilon^{1/2} \qquad \text{for all } t < \sigma \wedge \gamma_2(\varepsilon).$$

As a result, $\|X - x^*\|_{\sigma \wedge \gamma_2(\varepsilon)}^* \leq \varepsilon^{1/2} + 2\varepsilon \leq \gamma_1(\varepsilon)$. This shows (4.28) and completes the proof of Theorem 3.  $\square$

PROOF OF LEMMA 5.   By (4.18) and (4.20), the dynamics of $X$ on the intervals $I_k^1$ is given by

$$(4.32) \qquad X(t) = X(\eta_{k-1}) + W(t) - W(\eta_{k-1}) + (t - \eta_{k-1})(r + b),$$

where

$$(4.33) \qquad \begin{aligned} r_i := \ & \lambda_i - \sum_{j \in \mathcal{J}} \mu_{ij}\psi_{ij}^* - \alpha\delta_1 \sum_{j:(i,j)\in\mathcal{E}_p^-} \mu_{ij} \\ & + \delta_1 \sum_{j:(i,j)\in\mathcal{E}_p^+} \mu_{ij} + \delta_2 \sum_{j:(i,j)\in\mathcal{E}_p^c \cap \mathcal{E}_{\mathrm{ba}}} \mu_{ij} \end{aligned}$$

and

$$b_i = b_i(\mu, \beta) := - \sum_{j:(i,j)\in\mathcal{E}_p} \mu_{ij}\beta_{ij} - \sum_{j:(i,j)\in\mathcal{E}_p^c \cap \mathcal{E}_{\mathrm{ba}}} \mu_{ij}\beta_{ij}.$$

By (2.10), (4.10), (4.15)–(4.16) and (4.33), we have

$$(4.34) \qquad \sum_{i \in \mathcal{I}} r_i = \delta_2\Sigma_p^0 + \delta_1\Sigma_p^+ - \alpha\delta_1\Sigma_p^- \leq \frac{1}{2}(\delta_1\Sigma_p^+ - \alpha\delta_1\Sigma_p^-) < 0.$$

Note that by (4.6) and (4.1), $|\Delta X_e| \le 2\varepsilon$. From (4.20) and (4.21), and using (4.34), we thus obtain for $I_k^1$, $k \ge 1$,

$$-10\varepsilon \le X_e(\zeta_k) - X_e(\eta_{k-1})$$
$$\le W_e(\zeta_k) - W_e(\eta_{k-1}) + (e \cdot r + \|b\|)(\zeta_k - \eta_{k-1})$$
$$\le 2\varepsilon + (e \cdot r + c_1\varepsilon^2)(\zeta_k - \eta_{k-1}),$$

for $c_1 = \sum_{ij} \mu_{ij}$, and where we also used (4.6) and (4.17). Therefore, for $\varepsilon$ sufficiently small,

$$(4.35) \qquad |I_k^1| = \zeta_k - \eta_{k-1} \le m_1\varepsilon,$$

where $m_1 = 24/|e \cdot r|$. This proves part 1 of the lemma.

From (4.32) and (4.35), for $t \in I_k^1$ and sufficiently small $\varepsilon$,

$$\|X(t) - X(\eta_{k-1})\| \le 2\varepsilon + c_2(t - \eta_{k-1}) \le (2 + c_2 m_1)\varepsilon,$$

where $c_2 = 2\|r\|$. We therefore have

$$(4.36) \qquad \sup_{t \in [\eta_{k-1}, \zeta_k]} \|X(t) - X(\eta_{k-1})\| \le (2 + c_2 m_1)\varepsilon.$$

By (4.23), $\|X(t) - X(\zeta_k)\| \le 3\varepsilon$ for all $t \in I_k^2$, and taking into account a possible jump at $\eta_k$, we have

$$(4.37) \qquad \sup_{t \in [\zeta_k, \eta_k]} \|X(t) - X(\zeta_k)\| \le 5\varepsilon.$$

Since by (4.1) and (4.6), $\|X(0) - x^*\| \le \varepsilon$, part 2 of the lemma follows from (4.36) and (4.37).

In view of (4.22), (4.29) and (4.30), we have on $I_k^2$

$$(4.38) \qquad \begin{cases} Y(t) = 0, \\ Z(t) = -(X_e(t) - \nu_e - \theta_e)e_{j_0}, \\ \Psi_{ij}(t) = G_{ij}(X(t), \nu + \theta - Z(t)). \end{cases}$$

Define $\widetilde{X}(t) = X(t) - x^*$. From the definition of map $G$ (4.11)–(4.13), we have

$$(4.39) \qquad \begin{aligned} &G_{ij}(\widetilde{X}(t) + x^*, \nu + \theta - Z(t)) \\ &\quad = G_{ij}(x^*, \nu) + G_{ij}(\widetilde{X}(t) - \theta_e e_{i_0}, -Z(t)) + G_{ij}(\theta_e e_{i_0}, \theta). \end{aligned}$$

Due to (2.10), we have $G_{ij}(x^*, \nu) = \psi_{ij}^*$. Now consider the second term in (4.39). Using (4.38),

$$(4.40) \qquad \begin{aligned} G_{ij}(\widetilde{X}(t) - \theta_e e_{i_0}, -Z(t)) &= G_{ij}(\widetilde{X}(t) - \theta_e e_{i_0}, (X_e(t) - \nu_e - \theta_e)e_{j_0}) \\ &= G_{ij}(\widetilde{X}(t), \widetilde{X}_e(t)e_{j_0}) \\ &\quad + G_{ij}(-\theta_e e_{i_0}, -\theta_e e_{j_0}), \end{aligned}$$

where we used $\widetilde{X}_e = X_e - \nu_e$ due to $x_e^* = \nu_e$ (2.10). Finally, from (4.38)–(4.40), we have

$$
\begin{aligned}
(4.41) \qquad \Psi_{ij}(t) = \psi_{ij}^* &+ G_{ij}(\widetilde{X}(t), \widetilde{X}_e(t)e_{j_0}) \\
&+ G_{ij}(-\theta_e e_{i_0}, -\theta_e e_{j_0}) + G_{ij}(\theta_e e_{i_0}, \theta), \qquad t \in I_k^2.
\end{aligned}
$$

Define the map $H$ by

$$
(4.42) \qquad H_i(x) := -\sum_j \mu_{ij} G_{ij}(x, x_e e_{j_0}), \qquad x \in \mathbb{R}^I, \ i \in \mathcal{I},
$$

and the constant $H^\theta$ by

$$
(4.43) \quad H_i^\theta := -\sum_j \mu_{ij}[G_{ij}(-\theta_e e_{i_0}, -\theta_e e_{j_0}) + G_{ij}(\theta_e e_{i_0}, \theta)], \qquad i \in \mathcal{I}.
$$

By Assumption 1, $\sum_{j \in \mathcal{J}} \mu_{ij} \psi_{ij}^* = \lambda_i$. Hence, using (4.1), (4.38)–(4.43), we have

$$
\begin{aligned}
(4.44) \qquad \widetilde{X}(t) = \widetilde{X}(\zeta_k) &+ W(t) - W(\zeta_k) \\
&+ \int_{\zeta_k}^t H(\widetilde{X}(u))\,du + H^\theta(t - \zeta^k), \qquad t \in I_k^2,
\end{aligned}
$$

By (4.11)–(4.13) and (4.6), there exist constants $c_H > 0$ and $l_H > 0$, such that

$$
(4.45) \qquad \|H(x)\| \le c_H \|x\|, \qquad \|H^\theta\| \le l_H \varepsilon,
$$

for $\varepsilon$ sufficiently small. Therefore, applying (4.6), (4.45) and Lemma 5(2) to (4.44), we have from (4.23)

$$
3\varepsilon \le \|X(\eta_k) - X(\zeta_k)\| = \|\widetilde{X}(\eta_k) - \widetilde{X}(\zeta_k)\| \le 2\varepsilon + (c_H k m_2 + l_H)(\eta_k - \zeta_k)\varepsilon,
$$

for an appropriate constant $m_2$. The above implies

$$
\varepsilon \le (c_H k m_2 + l_H)(\eta_k - \zeta_k)\varepsilon.
$$

Therefore, for $k \ge 1$

$$
(4.46) \quad |I_k^2| = \eta_k - \zeta_k \ge \frac{1}{c_H k m_2 + l_H} \ge \frac{m_3}{k}, \qquad m_3 := \frac{1}{c_3 c_H m_2} < 1,
$$

where the constant $c_3$ satisfies $(c_3 - 1)c_H m_2 \ge l_H$, and where we assumed, without loss of generality, that $l_H > 1$ [cf. (4.45)]. This concludes the proof of Lemma 5. $\quad \square$

**5. Estimates on the probabilistic model.** In this section, we prove Theorem 1. We begin by introducing a rescaled version of the processes defined in Section 2 as follows. For $n \in \mathbb{N}$ and $t \geq 0$, let

$$(5.1) \qquad \bar{X}_i^n(t) = n^{-1} X_i^n(t), \qquad \bar{Y}_i^n(t) = n^{-1} Y_i^n(t), \qquad i \in \mathcal{I},$$

$$(5.2) \qquad \bar{Z}_j^n(t) = n^{-1} Z_j^n(t), \qquad \bar{\Psi}_{ij}^n(t) = n^{-1} \Psi_{ij}^n(t), \qquad i \in \mathcal{I}, \; j \in \mathcal{J}.$$

Denote $\bar{X}^n = (\bar{X}_i^n, i \in \mathcal{I})$, and use a similar convention for $\bar{Y}^n$, $\bar{Z}^n$ and $\bar{\Psi}^n$. Following a straightforward calculation, relations (2.1)–(2.5) can be rewritten in terms of the rescaled processes, as equations (5.3)–(5.7) below, holding for $n \in \mathbb{N}$ and $t \geq 0$:

$$(5.3) \qquad \bar{\Psi}_{ij}^n(t) = 0, \qquad (i,j) \in \mathcal{E}_{\mathrm{a}}^c,$$

$$(5.4) \qquad \begin{aligned} \bar{X}_i^n(t) = x_i^* + \bar{W}_i^n(t) + \lambda_i t \\ - \sum_{j \in \mathcal{J}} \mu_{ij} \int_0^t \bar{\Psi}_{ij}^n(s) \, ds, \qquad i \in \mathcal{I}, j \in \mathcal{J}, \end{aligned}$$

$$(5.5) \qquad \bar{Y}_i^n(t) + \sum_{j \in \mathcal{J}} \bar{\Psi}_{ij}^n(t) = \bar{X}_i^n(t), \qquad i \in \mathcal{I},$$

$$(5.6) \qquad \bar{Z}_j^n(t) + \sum_{i \in \mathcal{I}} \bar{\Psi}_{ij}^n(t) = \nu_j + \theta_j^n, \qquad j \in \mathcal{J},$$

$$(5.7) \qquad \begin{aligned} \bar{Y}_i(t) \geq 0, \qquad \bar{Z}_j(t) \geq 0, \\ \bar{\Psi}_{ij}^n(t) \geq 0, \qquad i \in \mathcal{I}, \; j \in \mathcal{J}, \end{aligned}$$

where we set

$$(5.8) \qquad \begin{aligned} \bar{W}_i^n(t) &:= n^{-1}[A_i^n(t) - \lambda_i^n t] \\ &- n^{-1} \sum_{j \in \mathcal{J}} \left[ S_{ij}^n \left( n \int_0^t \bar{\Psi}_{ij}^n(s) \, ds \right) - n \mu_{ij}^n \int_0^t \bar{\Psi}_{ij}^n(s) \, ds \right] \\ &+ (n^{-1} X_i^{0,n} - x_i^*) + (n^{-1} \lambda_i^n - \lambda_i) t - \sum_{j \in \mathcal{J}} (\mu_{ij}^n - \mu_{ij}) \int_0^t \bar{\Psi}_{ij}^n(s) \, ds \end{aligned}$$

and

$$(5.9) \qquad \theta_j^n = n^{-1} N_j^n - \nu_j.$$

The above equations resemble the dynamic fluid model studied in Section 4, and the proof of Theorem 1 will rely on the results of this section.

LEMMA 6. *Under any SCP, for any given $T \in (0, \infty)$, $\{n^{1/2} \|\bar{W}^n\|_T^*, n \in \mathbb{N}\}$ are tight random variables.*

PROOF. Relations (2.4), (2.5) and (2.7) imply that $0 \leq \bar{\Psi}_{ij}^n(t) \leq c_1$, where $c_1$ is a constant independent of $i, j, n$ and $t$. Hence by (2.15), the last three terms in (5.8) are bounded by $c_2(T+1)n^{-1/2}$, where $c_2$ is a constant independent of $i, j, n$ and $t$. Denote $\hat{A}_i^n(t) := n^{-1/2}(A_i^n(t) - \lambda_i^n t)$ and $\hat{S}_{ij}^n(t) := n^{-1/2}(S_{ij}^n(nt) - n\mu_{ij}^n t)$. Theorem 14.6 of [4] shows that $\{\hat{A}_i^n, n \in \mathbb{N}\}$ converges weakly to a Brownian motion (with zero mean and variance that depends on $i$), and that a similar statement holds for $\{\hat{S}_{ij}^n, n \in \mathbb{N}\}$. It follows that $\{|\hat{A}_i^n|_T^*, n \in \mathbb{N}\}$ and $\{|\hat{S}_{ij}^n|_{cT}^*, n \in \mathbb{N}\}$ are tight random variables, for each $i, j$, whenever $c$ is a constant that is independent of $n$. By (5.8), we obtain that

$$(5.10) \qquad n^{1/2}|\bar{W}_i^n|_T^* \leq |\hat{A}_i^n|_T^* + |\hat{S}_{ij}^n|_{c_1 T}^* + c_2(T+1).$$

As a result, $\{n^{1/2}|\bar{W}_i^n|_T^*, n \in \mathbb{N}\}$ are tight random variables, and the lemma follows. □

PROOF OF THEOREM 1. For $n \in \mathbb{N}$, let $\varepsilon^n = n^{-1/2}\log n$. By (2.15), for sufficiently large $n$,

$$(5.11) \qquad \|\theta^n\| \leq \varepsilon^n.$$

For $n \in \mathbb{N}$, let

$$(5.12) \qquad \tilde{\psi}_{ij}^n = \begin{cases} \psi_{ij}^* + \alpha\delta_1 + \beta_{ij}^n, & (i,j) \in \mathcal{E}_p^-, \\ \psi_{ij}^* - \delta_1 + \beta_{ij}^n, & (i,j) \in \mathcal{E}_p^+, \\ \psi_{ij}^* - \delta_2 + \beta_{ij}^n, & (i,j) \in \mathcal{E}_p^c \cap \mathcal{E}_{\text{ba}}, \\ 0, & (i,j) \in \mathcal{E}_p^c \cap \mathcal{E}_{\text{ba}}^c, \end{cases}$$

where $\beta_{ij}^n$ are constants chosen in such a way that for all sufficiently large $n$, one has $|\beta_{ij}^n|^2 \leq (\varepsilon^n)^2$, and $n\tilde{\psi}_{ij}^n$ has integer values, for each $i, j$ and $n$. Below, we write a system of equations for the processes $(\bar{X}^n, \bar{Y}^n, \bar{Z}^n, \bar{\Psi}^n)$ that uniquely defines them. We then let the processes $(X^n, Y^n, Z^n, \Psi^n)$ be defined through (5.1) and (5.2). These processes will then be shown to form a SCP, and to satisfy the statement of the theorem.

Fix $i_0 \in \mathcal{I}, j_0 \in \mathcal{J}$. Set $\eta_0^n = 0$, and consider the system of equations:

$$(5.13) \qquad \begin{cases} \bar{X}_i^n(t) = \bar{X}_i^n(\eta_{k-1}^n) + \bar{W}_i^n(t) - \bar{W}_i^n(\eta_{k-1}^n) + \lambda_i(t - \eta_{k-1}^n) \\ \qquad\qquad - \sum_{j \in \mathcal{J}} \mu_{ij} \int_{\eta_{k-1}^n}^t \bar{\Psi}_{ij}^n(s)\,ds, & i \in \mathcal{I}, t \in [\eta_{k-1}^n, \zeta_k^n), \\ \bar{\Psi}^n(t) = \tilde{\psi}^n, & t \in [\eta_{k-1}^n, \zeta_k^n), \\ \bar{Y}_i^n(t) = \bar{X}_i^n(t) - \sum_{j \in \mathcal{J}} \bar{\Psi}_{ij}^n(t), & i \in \mathcal{I}, t \in [\eta_{k-1}^n, \zeta_k^n), \\ \bar{Z}_j^n(t) = \nu_j + \theta_j^n - \sum_{i \in \mathcal{I}} \bar{\Psi}_{ij}^n(t), & j \in \mathcal{J}, t \in [\eta_{k-1}^n, \zeta_k^n), \end{cases}$$

where $\bar{W}^n$ is given by (5.8),

$$(5.14) \qquad \zeta_k^n = \inf\{t \geq \eta_{k-1}^n : \bar{X}_e^n(t) - \bar{X}_e^n(\eta_{k-1}^n) \leq -7\varepsilon^n\} \wedge \tau^n,$$

and

$$(5.15) \quad \begin{cases} \bar{X}_i^n(t) = \bar{X}_i^n(\zeta_k^n) + \bar{W}_i^n(t) - \bar{W}_i^n(\zeta_k^n) + \lambda_i(t - \zeta_k^n) \\ \qquad\quad - \sum_{j \in \mathcal{J}} \mu_{ij} \int_{\zeta_k^n}^t \bar{\Psi}_{ij}^n(s)\,ds, & i \in \mathcal{I}, t \in [\zeta_k^n, \eta_k^n), \\ \bar{Y}^n(t) = (\bar{X}_e^n(t) - \nu_e - \theta_e^n)^+ e_{i_0}, & t \in [\zeta_k^n, \eta_k^n), \\ \bar{Z}^n(t) = (\bar{X}_e^n(t) - \nu_e - \theta_e^n)^- e_{j_0}, & t \in [\zeta_k^n, \eta_k^n), \\ \bar{\Psi}^n(t) = G(\bar{X}^n(t) - \bar{Y}^n(t), \nu + \theta^n - \bar{Z}^n(t)), & t \in [\zeta_k^n, \eta_k^n), \end{cases}$$

where

$$(5.16) \qquad \eta_k^n = \inf\{t \geq \zeta_k^n : \|\bar{X}^n(t) - \bar{X}^n(\zeta_k^n)\| \geq 3\varepsilon^n\} \wedge \tau^n,$$

$$(5.17) \qquad \tau^n = \tilde{\tau}^n \wedge \sigma^n, \qquad \tilde{\tau}^n = \inf\{t \geq 0 : \|\bar{X}^n(t) - x^*\| \geq (\varepsilon^n)^{1/2}\},$$

$$(5.18) \qquad \sigma^n = \inf\{t \geq 0 : \|\bar{W}^n(t)\| \geq \varepsilon^n\},$$

and finally,

$$(5.19) \quad \begin{cases} \bar{X}_i^n(t) = \bar{X}_i^n(\tau^n) + n^{-1}(A_i^n(t) - A_i^n(\tau^n)), & t \geq \tau^n, \\ \bar{Y}_i^n(t) = \bar{X}_i^n(t), \qquad \bar{Z}_j^n(t) = n^{-1}N_j^n, \\ \quad \bar{\Psi}_{ij}^n(t) = 0, & i \in \mathcal{I},\ j \in \mathcal{J},\ t \geq \tau^n. \end{cases}$$

The above equations mimic the deterministic model (4.19)–(4.23). Note carefully that unlike in the case of (4.19)–(4.23), $\bar{W}^n$ should not be regarded here as data to the equations, since it depends on past values of $\bar{\Psi}^n$ (5.8). However, this will not create problem in applying the results of Section 4. Another difference is the definition of the processes for times $t \geq \tau^n$. Note that by definition of $\tau^n$, $\|\bar{W}^n(t)\| \leq \varepsilon^n$ holds for all $t < \tau^n$. This and (5.11) provide a bound that is similar to (4.6) over $[0, \tau^n)$. As a result, Lemma 4 implies that given the primitive processes $A_i^n$ and $S_{ij}^n$, the processes $\bar{X}^n, \bar{Y}^n, \bar{Z}^n$ and $\bar{\Psi}^n$ and the random times $\eta_k^n$, $\zeta_k^n$ and $\tau^n$ are uniquely defined by equations (5.13)–(5.19). It also follows from Lemma 4 that the processes $(\bar{X}^n, \bar{Y}^n, \bar{Z}^n, \bar{\Psi}^n)$ satisfy (5.3)–(5.7) over $[0, \tau^n)$. In turn, the processes $(X^n, Y^n, Z^n, \Psi^n)$ satisfy (2.1)–(2.5) on this interval. It is also easy to check that these equations are satisfied for $t \geq \tau^n$. We show that the processes $(X_i^n, Y_j^n, Z_j^n, \Psi_{ij}^n)$ take values in $\mathbb{Z}_+$. Since we have proved that (2.1)–(2.5) are satisfied, it suffices to show that $\Psi_{ij}^n$ take integer values. By construction of $\tilde{\psi}$, $\Psi_{ij}^n$ take integer values for $t \in [\eta_{k-1}^n, \zeta_k^n)$. On the intervals $[\zeta_k^n, \eta_k^n)$, by (4.22), $\Psi_{ij}^n$ will be the solution of the system of equations (4.11) with integer right-hand sides. In this case, a simple argument (the details of which are omitted) that uses the tree structure of $\mathcal{E}_{\mathrm{ba}}$ shows that $\Psi_{ij}^n$ are all integer valued.

To show that the constructed processes form a SCP, it remains to prove that for every $t$, $\Psi^n(t)$ is measurable on $\sigma\{X^n(s), A^n(s): s \leq t\}$ (cf. Definition 1). Fix $t$. We will show in steps (a)–(d) below that the value of $\Psi^n(t)$ is uniquely determined by the sample path $\Lambda[0, t] := \{X^n(s), A^n(s): s \in [0, t]\}$.

(a) By (2.2), the sample path $\Lambda[0, t]$ uniquely determines the sample paths $\sum_{j \in \mathcal{J}} S_{ij}^n(\int_0^\cdot \Psi_{ij}^n(u)\,du)$, $i \in \mathcal{I}$ on $[0, t]$.

(b) By (5.14), (5.16), $\Lambda[0, t]$ along with the value $\tau^n \wedge t$ uniquely determine the values $\eta_k^n \wedge t$, $\zeta_k^n \wedge t$, $k = 1, \ldots, K$. Thus, by (5.13), (5.15) and (5.19), $\Lambda[0, t]$ and $\tau^n \wedge t$ uniquely determine $\Psi^n$ on $[0, t]$. Equation (5.8), along with (a) above, shows that the same data, $\Lambda[0, t]$ and $\tau^n \wedge t$, uniquely determine $\bar{W}^n$ on $[0, t]$.

(c) We next show that $\Lambda[0, t]$ determines $\Psi^n$ and $\bar{W}^n$ on $[0, t)$. Let $(\Psi_1^n, \bar{W}_1^n)$, $(\Psi_2^n, \bar{W}_2^n)$ be two sample paths that correspond to the same data $\Lambda[0, t]$. Argue by contradiction and assume that on $[0, t)$ they do not agree. It follows from (b) that the corresponding values of $\tau^n$, that we denote by $\tau_1^n$ and $\tau_2^n$, do not agree, and that $\tau_1^n \wedge \tau_2^n < t$. Without loss of generality, assume that $\tau_1^n < \tau_2^n \wedge t$. Using (b) again, we have that $(\Psi_1^n, \bar{W}_1^n) = (\Psi_2^n, \bar{W}_2^n)$ on $[0, \tau_1^n]$. In particular,

$$\bar{W}_1^n(\tau_1^n) = \bar{W}_2^n(\tau_1^n). \tag{5.20}$$

Since $\widetilde{\tau}^n$ is defined in terms of $X^n$,

$$\widetilde{\tau}_1^n \wedge t = \widetilde{\tau}_2^n \wedge t. \tag{5.21}$$

Recall that $\sigma^n$ is the time when $\bar{W}^n$ leaves an open set of $\mathbb{R}^I$ (5.18). Thus, $\bar{W}_1^n(\tau_1^n)$ is either inside the open set, in which case $\tau_1^n = \widetilde{\tau}_1^n < \sigma_1^n \wedge \sigma_2^n$ and, therefore, by (5.21) $\tau_1^n = \tau_2^n$, or it is outside the open set, in which case $\tau_1^n = \sigma_1^n$, and by (5.20), we have that $\sigma_1^n = \sigma_2^n$ and $\tau_1^n = \tau_2^n$. In both cases, we obtain a contradiction to $\tau_1^n < \tau_2^n$. We conclude that statement (c) holds.

(d) By (5.8) and (a) above, $\bar{W}^n(t)$ is uniquely determined by $\Lambda[0, t]$ along with the values of $\Psi^n$ over $[0, t)$. Hence in view of (c), $\bar{W}^n(t)$ is determined by $\Lambda[0, t]$. Thus, the right continuity of $\bar{W}^n$ and the definition of $\sigma^n$ imply that $\sigma^n \wedge t$ is determined by $\Lambda[0, t]$. Hence, so is $\tau^n \wedge t$ and by (a), so is $\Psi^n(t)$.

Finally, we show that (2.16) and (2.17) hold for any fixed $T \in (0, \infty)$ and $\varrho > 1/2$. To this end, note that

$$\mathbb{P}(\sigma^n < T) \leq \mathbb{P}(\|\bar{W}^n\|_T^* \geq \varepsilon^n) = \mathbb{P}(n^{1/2}\|\bar{W}^n\|_T^* \geq \log n) \to 0,$$

by Lemma 6. On the event $\sigma^n \geq T$, we can use Theorem 3. On this event, for all $n$ sufficiently large, we have $T \leq \gamma_2(\varepsilon_n) \wedge \sigma^n$, thus Theorem 3 implies that $\int_0^T 1_{\{e \cdot Y^n(s) > 0\}}\,ds \leq \gamma_1(\varepsilon^n)$. Since $\varepsilon^n \to 0$ and $\gamma_1(0+) = 0$, (2.16) follows.

The second and third parts of Lemma 5 imply that there is a constant $C(T) < \infty$, independent of $n$, such that $\|\bar{X}^n - x^*\|_T^* \le C(T)\varepsilon^n$ on the event $\sigma^n \ge T$. On this event, we therefore have

$$n^{-\varrho}\|X^n - X^{0,n}\|_T^* \le C(T)n^{(1/2)-\varrho}\log n + n^{-\varrho}\|n^{-1}X^{0,n} - x^*\|,$$

where the last term on the above display converges to zero by (2.15). Since $\mathbb{P}(\sigma^n \ge T) \to 1$, (2.17) follows. This completes the proof of Theorem 1. □

**Acknowledgments.** We would like to thank the two referees for very useful comments.

DEPARTMENT OF ELECTRICAL ENGINEERING
TECHNION—ISRAEL INSTITUTE OF TECHNOLOGY
HAIFA 32000
ISRAEL
E-MAIL: atar@ee.technion.ac.il

DEPARTMENT OF MATHEMATICAL SCIENCES
CARNEGIE MELLON UNIVERSITY
PITTSBURGH, PENNSYLVANIA 15213
USA
E-MAIL: shaikhet@cmu.edu